\documentclass[final,epsfig,amsfont]{article}
\usepackage{latexsym}
\usepackage[dvips]{graphicx}
\newcommand{\eeq}{\end{equation}}
\newcommand{\beq}{\begin{equation}}
\newcommand{\nuq}[1]{\label{#1} \eeq}

\newtheorem{theorem}{Theorem}
\newtheorem{example}{Example}
\newtheorem{sequence}{Sequence}

\begin{document}
\title{Orthogonal polynomials of equilibrium measures supported on Cantor sets}
\author{
Giorgio Mantica \\
Center for Non-linear and Complex Systems,\\
Dipartimento di Scienze ed Alta Tecnologia, \\ Universit\`a dell'Insubria,\\
via Valleggio 11, 22100
Como, Italy. \\ Also at INFN sezione di Milano and CNISM unit\`a di Como.\\
{\small \em giorgio.mantica@uninsubria.it}
}
\date{}
\maketitle
\begin{abstract}
We study the orthogonal polynomials associated with the equilibrium measure, in logarithmic potential theory, living on the attractor of an Iterated Function System.
We construct sequences of discrete measures, that converge weakly to the equilibrium measure, and we compute their Jacobi matrices via standard procedures, suitably enhanced for the scope. Numerical estimates of the convergence rate to the limit Jacobi matrix are provided, that show stability and efficiency of the whole procedure. As a secondary result, we also compute Jacobi matrices of equilibrium measures on finite sets of intervals, and of balanced measures of Iterated Function Systems. 

These algorithms can reach large orders: we study the asymptotic behavior of the orthogonal polynomials and we show that they can be used to efficiently compute Green's functions and conformal mappings of interest in constructive function theory.

\end{abstract}


{\em MSC 2010 Class: 42C05 -- 31A15 -- 47B36 -- 81Q10 -- 30C30 }

\section{Introduction}
\label{sec1}

\subsection{Problem formulation and goals of this paper}
Orthogonal polynomials,  $\{p_j(\mu;s)\}_{j \in {\bf N}}$, of a positive Borel measure $\mu$ supported on a compact subset $E$ of the real axis are defined in a straightforward way
by the relation $\int p_j(\mu;s) p_m(\mu;s) d \mu(s) = \delta_{jm}$, where $\delta_{jm}$ is the Kronecker delta. The well known three-terms recurrence relation
\begin{equation}
\label{nor2}
   s p_j (\mu;s) = b_{j+1} p_{j+1}(\mu;s)
   + a_j p_j(\mu;s)  + b_{j} p_{j-1}(\mu;s),
\end{equation}
initialized by $b_0=0$ and
$p_{-1}(\mu;s) = 0$, $p_0(\mu;s) = 1$, can be formally encoded
in the Jacobi matrix $J(\mu)$:
\beq
   J(\mu) :=
           \left(   \begin{array}{ccccc}  a_0  & b_1 &      &      &       \cr
                               b_1  & a_1 & b_2  &      &       \cr
                      &   \ddots    &    \ddots      & \ddots &\cr
                               \end{array} \right) .
 \nuq{jame}
For compact support $E$ the moment problem is determined \cite{ach}, and the matrix $J(\mu)$
is in one--to--one relation with the measure $\mu$.

While originally introduced for applications (quadratures, optimal control) the r\^ole of orthogonal polynomials in harmonic analysis, analytical functions and potential theory soon emerged \cite{akhie1,gero}, and appear clearly in their asymptotic properties for large order, beautifully described in the by--now classical book \cite{stahl}. These relations are particularly intriguing in the case of measures supported on Cantor sets, the object of this paper. Consider in fact limits
such as the {\em ratio asymptotics} $p_{j+1}(\mu;z)/p_j(\mu;z)$, or the {\em j-th root asymptotics}, $|p_j(\mu;z)|^{1/j}$, where $z$ is a point in the complex plane and the order $j$ tends to infinity. Under well specified conditions \cite{stahl}, these limits exist and yield the Green's function, $g(E;z)$, of the Dirichlet problem for the complement of the set $E$ \cite{ed,ran0}. In turn, the Green's function is related to an additional measure, that appears in two different forms, $\nu_\mu$ and $\nu_E$, that may, or may not, exist and coincide.

On the one hand, $\nu_\mu$ is the {\em counting measure} of the zeros of
the orthogonal polynomials: letting $\xi^j_l$, for $l=1,\ldots,j$, be the zeros of $p_j(\mu;z)$ and $D_x$ be the unit mass, atomic (Dirac) measure located at the point $x$, the measure $\nu_\mu$ is defined by
 \begin{equation}
 \label{pote0}
    \nu_\mu = \lim_{j \rightarrow \infty} \frac{1}{j}  \sum_{l=1}^j D_{\xi^j_l},
 \end{equation}
where convergence is meant in the weak * sense.

On the other hand, if $E = \mbox{supp} (\mu)$ is a compact subset of the complex plane $\bf C$, $\nu_E$ is the {\em electrostatic equilibrium measure} for a charge distributed on $E$, with a logarithmic law of repulsion. In fact, let
$\sigma$ any Borel probability measure, also supported on $E$. The potential $V(\sigma;z)$, generated by $\sigma$ at the point $z$ in $\bf C$, is
 \begin{equation}
 \label{pote1}
    V(\sigma;z) := - \int_E \log |z-s| \; d \sigma(s).
 \end{equation}
The electrostatic energy ${\cal E}(\sigma)$ of the distribution $\sigma$ is given by the integral of $V(\sigma;z)$:
 \begin{equation}
 \label{pote2}
  {\cal E}(\sigma) := \int_E  V(\sigma;u) \; d\sigma(u)  = - \int_E \int_E \log |u-s| \; d\sigma(s)d\sigma(u).
 \end{equation}
The {\em equilibrium measure} $\nu_E$ associated with the compact domain $E$ is the unique measure that minimizes the energy ${\cal E}(\sigma)$, when this latter is not identically infinite \cite{ran0,ed}. In this case, $\mbox{Cap}(E) := e^{-{\cal E}(\nu_E)}$  defines the capacity of the set $E$, and
the Green's function can be written as
\begin{equation}
 \label{pote8}
  g(E;z) = - V(\nu_E;z) - \log (\mbox{Cap}(E)) = - V(\nu_E;z) + {\cal E}(\nu_E).
 \end{equation}

Observe that in the above definition $\nu_E$ depends only on the set $E$. It coincides with $\nu_\mu$ in the so called {\em regular} case:  measures $\mu$ not too thin on any part of their support are regular---see  \cite{stahl} for the exact definition. The measures studied herein will all be regular, so to enable us to use both characterizations, $\nu_\mu = \nu_E$, and the existence of root asymptotics.

Our goal in this paper is to study the measure $\nu_E$ and its orthogonal polynomials $p_j(\nu_E;z)$, by devising a reliable computational scheme for the associated Jacobi matrix, $J(\nu_E)$,
a task which, to the best of our knowledge, is still to be performed
when $E$ is a Cantor set generated by Iterated Function Systems (IFS) \cite{papmor,hut,dem,ba2,recur}, to be described in the following.
Furthermore, we want to analyze the convergence rate in the root asymptotics defined above.
The computation of the Jacobi matrix is a fundamental problem in numerical analysis, for a wealth of reasons \cite{gautschi,tom}. In this paper we will try to add a further one:
when combined with root asymptotics, it yields efficient algorithms to compute the Green's function $g(E;z)$, the electrostatic potential $V(\nu_E;z)$ and conformal mappings of interest in constructive function theory.

\subsection{Background and motivations}
While many problems in classical orthogonal polynomials have found a complete solution, both from the analytical and the computational viewpoint, much is still to be discovered for non-classical orthogonal polynomials supported on Cantor sets on the real line\footnote{Since the terminology {\em semi-classical orthogonal polynomials} \cite{maroni} also exists, it is tempting to call {\em quantum orthogonal polynomials} those associated with this case. This would not be fully inappropriate, because quantum mechanical phenomena originated by these measures 
have been termed {\em quantum intermittency}
\cite{physd1,physd2,etna,poin1,poin2,serguey}.}, despite considerable progress has been made in the last thirty years or so, in part for applications to quantum mechanics.
In fact, the Jacobi matrix
$J(\mu)$ can be seen as an operator acting on $l^2({\bf Z}_+)$, the space of square summable sequences. For instance, choosing $b_j=1/2$ for all $j$ and writing formally $a_j = F(j)$, with $F$ a potential function, yields a {\em discrete Schr\"{o}dinger operator} \cite{cycon,testard,last,last2}. In so doing, $\mu$ becomes the {\em spectral measure} of $J(\mu)$ associated with the first basis vector of $l^2({\bf Z}_+)$, while the equilibrium measure $\nu_E$ is known as the {\em density of states}, a measure which plays a fundamental r\^ole in many physical properties of the system \cite{ldos,reedsim}.

Quite naturally, the question arises on the relations between the properties of the two sequences $\{a_j\}_{j \in {\bf N}}$, $\{b_j\}_{j \in {\bf N}}$, and those of $\mu$. For instance, much is known about measures in the Nevai class $N(a,b)$ ({\em i.e.} those for which $\{a_j\}_{j \in {\bf N}}$ and $\{b_j\}_{j \in {\bf N}}$ tend to finite limits $a$ and $b$), as well for Jacobi matrices that are asymptotically periodic \cite{jeff1,vanas,lubin}. Recently, the link between these limits and the classical Sz\"ego asymptotics of orthogonal polynomials has been fully clarified \cite{damani,simonratio}.
Next in complexity comes the case of discrete Schr\"odinger operators with {\em almost periodic} potentials \cite{last,testard,baker,belli,mijeff,russell,physd1,physd2}. Typically, it is found that the associated spectral measures are singular continuous, and in certain cases (like the so--called Fibonacci matrix \cite{fibo1,fibo2}) supported on Cantor sets with self--similar geometry.
Conversely, one can start with IFS {\em balanced} (not equilibrium) measures on such sets, and ask what are the properties of the associated Jacobi matrices---that can be computed numerically \cite{hans2,cap,arxiv}.
It is still an open problem to assess whether almost periodicity of some sort characterizes these matrices, as conjectured in \cite{physd1}. Results in this direction can be found in the theory of equilibrium measures on finitely many intervals (see eg. \cite{widom,nuttall,alphonse,sasha1,franz,franzimrn,franz0,chris,chris2}).
Presumably the proof is to be found in the properties of the equilibrium measure $\nu_E$, and the study of $J(\nu_E)$ might be a good starting point.
Numerical investigations have indeed served as seeds of serendipitous discovery in various problems of harmonic analysis on fractals \cite{str0,str0b,stric,stric0}.

Finally, equilibrium measures are studied in constructive function theory. For instance, it has been shown that the behavior of $V(\nu_E;z)$ on the complement of a compact set $E$ depends on the smoothness of $\nu_E$ \cite{totik,andriev,andriev3}. In these investigations, Cantor sets play a special r\^ole, since they provide nice applications of the abstract theorems. It is then obvious that numerical experiments can serve both as an illustration of these theorems and as a stimulus of new results.

\subsection{Outline of the paper and summary of results}

The typical construction of Cantor sets $E$ on the real line via IFS
is reviewed in Sect. \ref{sec-ifs}, where we introduce two families of IFS: the first is composed of affine maps, like those generating the middle--third Cantor set (Example 1 for the numerical tests of this paper), the second consists of non--linear IFS yielding Julia sets (Example 2). This construction produces a sequence of finite unions of intervals, $\{E_n\}_{n \in {\bf N}}$, converging to $E$ when $n$ tends to infinity. IFS balanced measures (not to be confused with equilibrium measures) on $E$ are also defined in Sect. \ref{sec-ifs}.

Each $E_n$ is a compact set, that carries a unique equilibrium measure, $\nu_n$. Numerical approximations of $\nu_n$, as the weak limit of a sequence of {\em discrete} measures, composed of a finite number of atoms, can be produced by a judicious use of Gaussian integration and the solution of a system of non--linear equations \cite{dolo}. This is described in Sect. \ref{sec-equil} and in Appendix A.

Following suit, one needs algorithms to compute the Jacobi matrix associated with a finite number of atoms: in Sect. \ref{sec-gragg}, we test four different known techniques to this scope. Our benchmark is a sequence of discrete measures converging to the IFS balanced measure on a Cantor set. As a result of this experiment, we select the method RKPW introduced by Gragg and Harrod \cite{gragg} as the best for our purpose. We bring minor modifications to its standard version, described in Appendix B, that reduce its computational complexity, to make it affordable also when considering large numbers of atoms. We prove that the modified algorithm can also be seen as a technique to add a finite number of atomic measures to the Jacobi matrix of {\em any} arbitrary measure (as Fisher's method \cite{hans3}, but in a stable fashion).

This procedure is applied to the computation of the Jacobi matrix of a balanced IFS measure in Sect. \ref{sec-jacbala}, and of the equilibrium measure $\nu_n$ on $E_n$ in Sect. \ref{sec-jeq}.
We introduce a fundamental quantity, $N_\epsilon$, as the rank of the largest truncated Jacobi matrix that is computed with absolute component-wise error less than $\epsilon$. This quantity helps us to control both the maximum precision attainable by the algorithm and its computational complexity.

The sequence of measures $\{\nu_n\}_{n \in {\bf N}}$ is
convergent, in the weak * topology, to $\nu_E$, the equilibrium measure on $E$. Therefore, in Sect. \ref{sec-jeqattra}, we tackle the limit process of letting the order $n$ go to infinity, to compute $J(\nu_E)$, the Jacobi matrix of the equilibrium measure on a Cantor set $E$. We describe in detail the algorithms and we present numerical results for the two examples mentioned before---ternary Cantor and Julia sets. We study the convergence properties as a function of the precision required, the order $n$ of the IFS construction and the number of Gaussian points employed. The success of our numerical technique lies in the fact that convergence takes place orderly, from the top entries of the Jacobi matrix downwards, extending to very large indices, due to the slow growth of the numerical error.

The rest of the paper outlines two applications of the Jacobi matrices $J(\nu_n)$ and $J(\nu_E)$ so constructed: root asymptotics and conformal mappings. First, in Sect. \ref{sec-capa}, we consider the asymptotic behavior of the sequence of orthogonal polynomials $\{p_j(\nu_E;z)\}_{j \in {\bf N}}$. We focus on $j$-th root asymptotics, that yields in the limit the {\em real} Green's function of the electrostatic problem for the set $E$. The usual concept can be extended, following \cite{widom,thouless,johnson}, to a {\em complex} Green's function, whose real part is the Lyapunov exponent, and whose (harmonic conjugate) imaginary part extends to the complex plane the rotation number for discrete Sturm--Liouville operators \cite{souil,johnson}. Regular root asymptotics is proven to hold also in this meaning, in the cases studied; in addition, we perform numerical experiments that measure the rate of convergence to the asymptotic limit. As a consequence of this investigation, we demonstrate that root asymptotics is an efficient numerical tool to compute the complex Green's function. Finally, in Sect. \ref{sec-confor}, we describe an algorithm for a
conformal mapping of the external of the Cantor set $E$, to the external of the unit disk \cite{andriev}, that employs the Jacobi matrices derived in the preceding sections. The conclusions briefly mention further examples where these techniques can be profitably applied.

\section{IFS attractors and Balanced Measures}
\label{sec-ifs}
Iterated Function Systems (IFS) \cite{papmor,hut,dem,ba2,recur,palle} provide
a convenient construction of Cantor sets. In the simplest setting,
they are collections of contractive maps $\phi_m : {\bf R} \rightarrow
{\bf R}$, $m = 1, \ldots, M$: for any $m$ there exists $\delta_m < 1$, such that
$|\phi_m(s)-\phi_m(t)| \leq \delta_m |s-t|$. There exists a unique
set $E$, called the {\em attractor} of the IFS, that solves
the equation
 \begin{equation}
 \label{attra}
    E=\bigcup_{m=1,\ldots ,M}\;\phi_m(E) :=  \Phi (E).
 \end{equation}
In the above, we have also defined the operator $\Phi$, acting on the set of compact subsets of ${\bf R}$. This space is complete in the Hausdorff metric, and
$\Phi$ is there contractive. Therefore, the attractor $E$ is also the limit of the sequence $\Phi^{n}(E_{0})$, where $E_{0}$ is any non-empty compact set:
$
E = \lim_{n \to \infty} \Phi^{n}(E_{0}).
$

The set of measures supported on $E$ is rich and wide. We will restrict our consideration to two kinds of measures of mathematical and physical significance: {\em equilibrium measures} and {\em balanced measures}. The former having been defined in the Introduction, let us now describe the latter.
They are obtained associating a probability, or weight, $\pi_m>0$, $m = 1, \ldots, M$, to each IFS map: $\sum_m \pi_m = 1$. For any such choice of weights, a unique positive measure $\mu$ on $E$ satisfies the equation
\begin{equation}
\label{bala}
   \int f \; d\mu \; =
   \sum_{m=1}^{M}
   \; \pi_{m}
   \; \int   \;
 (f \circ \phi_{m}) \; d\mu ,
\end{equation}
for any continuous function $f$. For instance, consider
the set of one--dimensional affine maps of the form:
\begin{equation}
\label{mappi}
    \phi_m (s) = \delta_m (s - \gamma_m) + \gamma_m,  \;\;  m = 1, \ldots, M ,
\end{equation}
where $\delta_m$ are real numbers between zero and one, called {\em contraction ratios}, and $\gamma_{m}$ are real constants, the fixed points of the maps.
Under these conditions, the attractor $E$ is a finite or infinite collection of intervals, or a Cantor set  \cite{nalgo2,intjns}.
\begin{example} {\bf Ternary Cantor Set}. Let $M=2$, $\delta_1=\delta_2=1/3$, $\gamma_1=0$, $\gamma_2=1$,
$\pi_1=\pi_2=1/2$. The attractor of this IFS is the middle--third Cantor set and the balanced measure $\mu$ is the Devil's staircase measure.
\label{examp-cantor}
\end{example}
Equation (\ref{bala}) can be rewritten as $\int f d \mu = \int (T f) d \mu$, where we have introduced the transfer operator $T$. By going to the dual space of Borel probability measures, ${\cal M}$, this equation is equivalent to
$
   T^* \mu = \mu,
$
where $T^*$ is the adjoint operator of $T$, also known as the Perron--Frobenius operator. This operator is contractive in the Hutchinson--Wesserstein metrics, so that the sequence of discrete measures $\mu_n := (T^*)^n\mu_0$, for any $\mu_0$, converges in the complete space ${\cal M}$. This permits to define a first sequence of discrete measures. As above,
let $D_x$ be the Dirac measure located at the point $x$.
\begin{sequence}
\label{ex-1}
{\bf Discrete measures converging to the balanced measure of an IFS.}
Let $\mu_0 = D_{\frac{1}{2}}$. For any $n \in {\bf N}$, $\mu_n:= (T^*)^n \mu_0$ is a discrete measure composed of $M^n$ atoms:
\begin{equation}
\label{bala3}
  \mu_n = \sum_{i=1}^{M^n} w_i D_{x_i},
\end{equation}
where $x_i$ and $w_i$ are easily constructed by recurrence, via eqs. (\ref{bala}),(\ref{mappi}). The sequence $\{\mu_n\}_{n \in {\bf N}}$ converges to $\mu$, the balanced measure on the IFS attractor.
\end{sequence}

A second interesting family of I.F.S. maps consists of the real inverse roots of polynomials in a complex variable.
\begin{example}
{\bf Real Julia set of a quadratic mapping}. Consider the IFS composed of $M=2$ non--linear maps
\beq
  \phi_\pm (x) = \pm \sqrt{x + \lambda},
  \label{eq-julia1}
  \eeq
where $\lambda \geq 2$ is a real parameter, with associated probabilities $\pi_\pm = 1/2$.
\label{exa-julia}
\end{example}

Observe that the two maps in eq. (\ref{eq-julia1}) are the inverse branches of the quadratic transformation $z \rightarrow z^2 - \lambda$, so that the IFS attractor $E$ is the Julia set of this map.
When $\lambda = 2$ this set is the interval $[-2,2]$ and the balanced measure is the Chebyshev measure, {\em i.e.} the equilibrium measure on $E$: $\mu = \nu_E$. This remarkable coincidence is general:
\begin{theorem}
\label{th-julia}
For any $\lambda \geq 2$ the Julia set of $z \rightarrow z^2 - \lambda$ is a subset of the real line, whose equilibrium measure $\nu_E$ coincides with the balanced measure $\mu$ of the IFS in Example \ref{exa-julia}.
\end{theorem}
For proof and theory of Julia sets see \cite{jul,fat,brol,bessi}.
Out of this beautiful theory we need to recall a second fact: because of the renormalization relation $p_{2j}(\nu_E;\phi_\pm (x)) = p_j (\nu_E;x)$, the Jacobi matrix $J(\nu_E)$ can be computed via simple recursions \cite{belli,mijeff}, $a_j=0$ for all $j$, while $b_1=\lambda$,
$b_{2j} = b_{j}/b_{2j-1}$ and $b_{2j+1}=\lambda-b_{2j}$. We will use these relations in Sect. \ref{sec-jeqattra}.

\section{Equilibrium measures on IFS Attractors}
\label{sec-equil}

The logarithmic capacity of the attractor $E$ of an IFS has been studied in \cite{ran1,ran,dolo}. We now compute numerically its equilibrium measure, following \cite{dolo}.
Let $E_0$ be the convex hull of $E$, that can be easily be identified as the interval
 $   E_0 =  [\gamma_1,\gamma_M]$,
where we have ordered the IFS maps according to increasing values of their fixed points:  $\gamma_j < \gamma_{j+1}$, for any $j=1,\ldots,M-1$.
The set $E$ is the limit, in the Hausdorff metric, of the sequence of compact sets $E_n:= \Phi^n(E_0)$, see eq. (\ref{attra}), each of which the union of $N$ disjoint, closed intervals $E_n^i$:
\begin{equation}
\label{eq-int1}
    E_n = \Phi^n(E_0) = \bigcup_{i=1}^{N} E_n^i.
\end{equation}
The maximum cardinality, $N=M^n$, is met in the case of {\em fully disconnected} IFS, ({\em i.e.} those for which the intervals $\{\phi_m(E_0)\}_{m=1}^{M}$ are pairwise disjoint).  We denote the intervals in eq. (\ref{eq-int1}) as $E_n^i := [\alpha_i,\beta_i]$, dropping for simplicity the generation index $n$ when confusion cannot occur. The idea of approximating $E$ by the nested sequence of sets $E_n$ has proven to be useful also theoretically \cite{andriev}: thanks to basic properties of the logarithmic potential, capacity and Green's function can be continuously obtained in the limit $n \rightarrow \infty$. The same holds for equilibrium measures and Jacobi matrices, as we will see momentarily.

The equilibrium problem for a finite union of $N$ intervals $[\alpha_i,\beta_i]$, $i = 1,\ldots,N$, has a well known analytical solution
\cite{akhie2,widom,nuttall,alphonse,sasha1,franz0,franz,franzimrn,chris}: define the polynomial $Y(z)$,
 \begin{equation}
 \label{meas2}
    Y(z) = \prod_{i=1}^{N} (z-\alpha_i)(z-\beta_i),
 \end{equation}
and its square root, $\sqrt{Y(z)}$, as the one which takes real values for $z$ real and large. Also, let the real numbers $\xi_i$ belong to the open intervals $(\beta_i,\alpha_{i+1})$, for $i=1,\ldots,N-1$. Define $Z(\xi;z)$ as the monic polynomial of degree $N-1$ with roots at all $\xi_i$'s:
 \begin{equation}
 \label{meas3}
    Z(\xi;z) =  \prod_{i=1}^{N-1} (z-\xi_i).
 \end{equation}
Then, there exists a unique set of values $\{\zeta_i\}_{i = 1,\ldots,N-1}$ that solve the set of coupled, non--linear equations
 \begin{equation}
 \label{meas5}
    \int_{\beta_i}^{\alpha_{i+1}} \frac{Z(\zeta;s)}{\sqrt{|Y(s)|}} \; ds = 0, \;\; i = 1,\ldots,N-1.
 \end{equation}
A stable technique for the solution of the non-linear equations (\ref{meas5}) has been presented in \cite{dolo}. To reduce its memory and computer time requirements we developed a different technique, that is described in Appendix A, not to interrupt here the natural flow of arguments.

When the solution $\{\zeta_i\}_{i=1,\ldots,N-1}$ of equations (\ref{meas5}) is known, the equilibrium measure on $E_n$, denoted by $\nu_n := \nu_{E_n}$, can be computed as the absolutely continuous measure (with respect to the Lebesque measure on $E_n$), consisting of the sum of $N$ measures supported on each of the $N$ intervals $E_n^i$:
 \begin{equation}
 \label{meas6}
     d \nu_n (s) = \frac{1}{\pi} \sum_{i=1}^N \chi_{[\alpha_i,\beta_i]} (s)\;  \frac{|Z(\zeta;s)|}{\sqrt{|Y(s)|}} \; ds.
 \end{equation}
The properties of the orthogonal polynomials of $\nu_n$, and of the Green's function of $E_n$ have been studied in \cite{akhie2,widom,nuttall,sasha1,franz0,franz,jeff2}. In particular, the algebraic approach of \cite{franzimrn, alphonse,alphonse2} can be turned into a symbolic computation of these quantities. In this paper, to the contrary, we adopt a numerical strategy.

In fact, when considering the $i$-th interval composing $E_n$, the function $Y(s)$ can be factored as $Y(s) = (s-\alpha_i)(\beta_i-s) \tilde{Y_i}(s)$, with obvious meaning of the function $\tilde{Y}_i(s)$. Therefore, the measure with density $\frac{|Z(\zeta;s)|}{\pi \sqrt{|Y(s)|}}$ with respect to the Lebesgue measure on $[\alpha_i,\beta_i]$, can also be seen as the absolutely continuous measure, with respect to the Chebyshev measure on the same interval, with density $\frac{|Z(\zeta;s)|}{\sqrt{|\tilde{Y}_i(s)|}}$:
 \begin{equation}
 \label{meas6t}
       \frac{|Z(\zeta;s)|}{\sqrt{|Y(s)|}} \; ds = \frac{|Z(\zeta;s)|}{\sqrt{|\tilde{Y}_i(s)|}} \frac{d s}{\sqrt{(s-\alpha_i)(\beta_i-s)}}, \;\; i = 1,\ldots,N.
 \end{equation}
This fact enables us to define a second sequence of point measures:
\begin{sequence}
{\bf Discrete measures converging to the equilibrium measure}.
Let $\theta^G_l$, $l=1,\ldots,G$, be the Gaussian points of order $G$ for the Chebyshev measure on $[-1,1]$, and let $\psi_n^i$ the affine map that takes $[-1,1]$ unto $E_n^i$, for $i=1,\ldots,N$. For any $n$ and $G \in \bf{N}$, define the discrete measure $\nu_n^G$ by:
\begin{equation}
 \label{eq-meas20}
     \nu_n^G  =
     \frac{1}{G}
      \sum_{i=1}^N \sum_{l=1}^G
     \frac{|Z(\zeta;\psi_i(\theta^G_l))|}{\sqrt{|\tilde{Y_i}(\psi_n^i(\theta^G_l))|}} \;
     D_{\psi_n^i(\theta^G_l)}.
 \end{equation}
When $G$ tends to infinity, $\nu_n^G$ tends weakly to $\nu_n$, the equilibrium measure on $E_n$.
\end{sequence}

\begin{figure}
\centerline{\includegraphics[width=.6\textwidth, angle = -90]{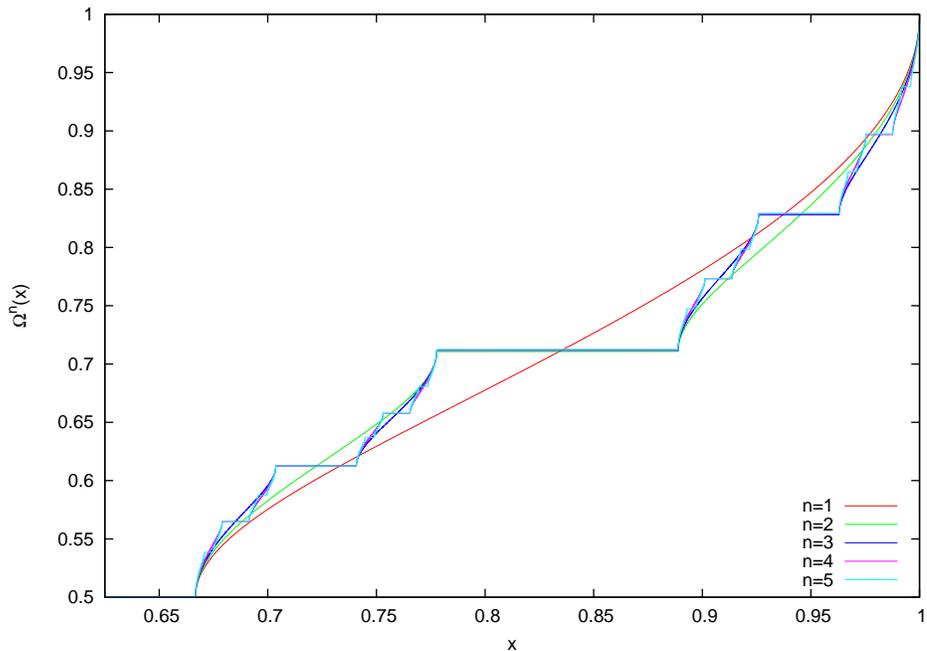}}
\caption{Equilibrium measures $\Omega_n(x)$ versus $x$ at various generation numbers $n$, for the IFS of Example \ref{examp-cantor}. Because of symmetry, only half of the picture is shown.}
\label{erra2}
\end{figure}

In Figure \ref{erra2} we plot $\Omega_n^G(x) := \int_{\alpha_1}^x d \nu_n^G$, the distribution function of the discrete measure $\nu_n^G$, versus $x$, for $n=1,\ldots,5$, when $G$ is kept constant and large, in the case of Example \ref{examp-cantor}.
Since $\nu_n$ tends to the equilibrium measure $\nu_E$ on $E$ when $n$ tends to infinity, it is clear that letting both $n$ and $G$ grow simultaneously\footnote{One can also choose a different number of Gaussian points for each set $E_n^i$. This is useful when their lengths are vastly different, like in affine IFS with different contraction ratios, or in the Julia set case.} we can devise a sequence of discrete measures $\nu_n^G$ converging to $\nu_E$. An efficient procedure to rule their combined increase can be derived by the convergence of the corresponding Jacobi matrix entries, as we will show in Sect. \ref{sec-jeqattra}.

Finally, observe that the distribution of atoms composing $\nu^G_n$ is {\em not} a proper Gaussian measure associated with $\nu_n$ (one that is obtained by diagonalization of a truncation of the Jacobi matrix---see below, in Sect. \ref{sec-gragg}); rather, it is a more easily computable approximation, in the same spirit of near--optimal distributions of Riesz energy points on manifolds \cite{doug}.

\section{Algorithms to compute the Jacobi matrix of discrete measures}
\label{sec-gragg}

In the previous sections we have described two sequences of discrete measures, $\{\mu_n\}_{n \in {\bf N}}$ (eq. \ref{bala3}) and $\{\nu_n^G\}_{G \in {\bf N}}$ (eq. \ref{eq-meas20}) that converge to $\mu$ and $\nu_n$, respectively. Algorithms to compute the Jacobi matrix of discrete measures number in the many. Yet, when applied to sets of atoms of large cardinality, whose locations approach a Cantor set, they may suffer from numerical instabilities.
Therefore, in this section we test four algorithms, in search for the best to be applied to our problem. They are: {\em a:} an implementation of the Stieltjes/Lanczos technique, where integration is performed by a finite summation; {\em b:} Fischer's approach \cite{hans3} of adding an atomic measure to a second measure (known via its Jacobi matrix); {\em c:} Laurie's quotient-difference method \cite{dirk9}; {\em d:} Gragg and Harrod's algorithm RKPW \cite{gragg}, based on plane rotations.

{\em Numerical experiment 1. Consider the weakly convergent sequence $\{\mu_n\}_{n \in {\bf N}}$, Sequence \ref{ex-1}, in the case of Example \ref{examp-cantor}. For each $n$, compute the associated Jacobi matrix $J(\mu_n)$, of size $N = M^n$, (here, $M=2$) via each of the four algorithms above. Next, diagonalize this latter:
\beq
  J(\mu_n) \; {\bf u}_i = \lambda_i {\bf u}_i,
  \label{eq-goluwel1}
  \eeq
to obtain, via the Golub--Welsh algorithm \cite{golw}, eigenvalues $\lambda_i$ and squared first components of the normalized eigenvectors, $v_i = (e_0,{\bf u}_i)^2$, for $i=1,\ldots,N$}.

\begin{figure}
\centerline{\includegraphics[width=.6\textwidth, angle = -90]{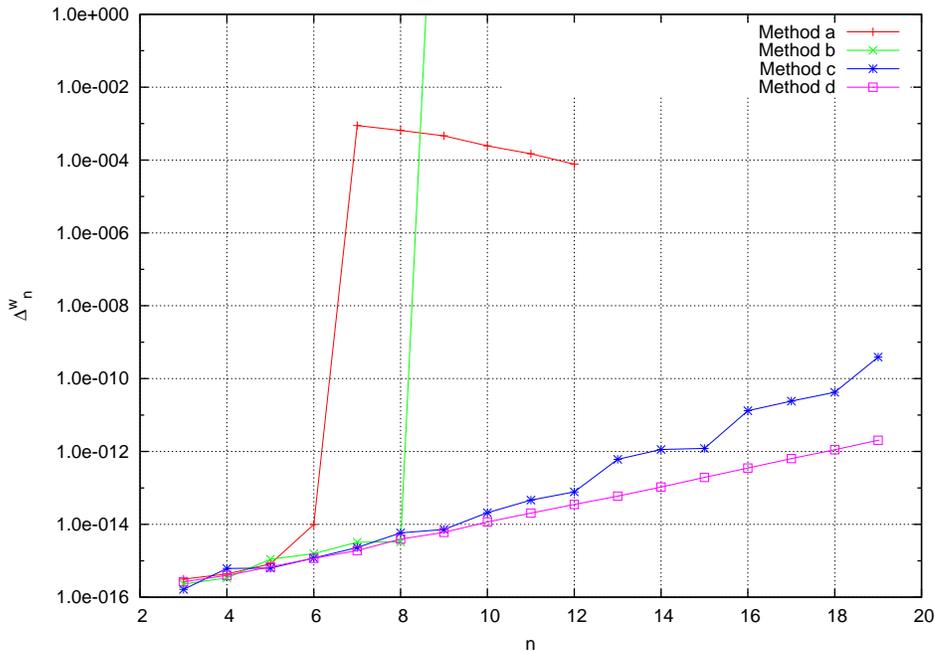}}
\caption{$L^1$ errors $\Delta^w_n$ in the reconstruction of the weights of the discrete measure $\mu_n$, versus generation level $n=3,\ldots,19$, for the IFS of Example \ref{examp-cantor}. Methods a) to d) are described in the text.}
\label{fig-err1}
\end{figure}

The last set of values should exactly reproduce the points and weights, $x_i$ and $w_i$, respectively, defining the discrete measure $\mu_n$. The $L^1$ errors in this reconstruction are $\Delta^x_n := 1/N \sum_i |x_i -\lambda_i|$ and $\Delta^w_n := 1/N \sum_i |w_i -v_i|$. In Figure \ref{fig-err1} we plot $\Delta^w_n$ versus $n$, for the four methods listed above.  The errors $\Delta^x_n$ are smaller than $\Delta^w_n$ and are not plotted.
We observe that the Lanczos technique gives fully unreliable results at $n=7$, while Fischer's method breaks down (negativity of a positive definite quantity) at $n=9$. To the contrary, methods {\em c} and {\em d} never break down in the range explored. Method {\em d} seems to be the best, featuring an error that increases as a power--law in $N$, with exponent smaller than one. Since the results of this experiment are typical, we have adopted Gragg and Harrod's algorithm (RKPW henceforth) in the computations of this paper. In Appendix B we show that
it can be conveniently programmed on a parallel computer, and that its computational complexity can be significantly reduced, when considering finite Jacobi matrices of rank much less than the number of atoms in the associated measure.

\section{The Jacobi matrix of balanced IFS measures}
\label{sec-jacbala}

In the previous section we have computed
the Jacobi matrices of the measures $\{\mu_n\}_{n \in {\bf N}}$, that converge to the balanced measure $\mu$ of an IFS. This provides a new algorithm for the difficult numerical problem of computing the Jacobi matrix $J(\mu)$. Since this latter can also be obtained by different, non iterative techniques \cite{hans2,cap,arxiv}, we can measure the convergence of the matrix entries of $J(\mu_n)$ as a function of the order $n$, to assess the relative performance the new algorithm. For this, we further analyze the results of Experiment 1.

\begin{figure}
\centerline{\includegraphics[width=.6\textwidth, angle = -90]{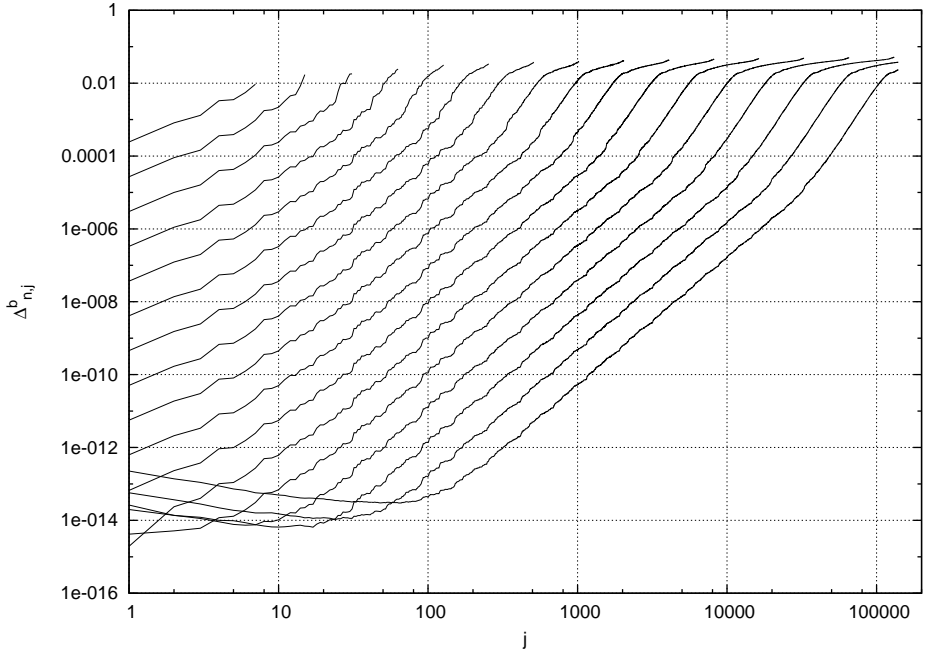}}
\caption{Differences $\Delta^b_{n,j}$ in the computation of the Jacobi matrix of the IFS in Example \ref{examp-cantor}, plotted versus $j$.
Data with the same $n$ are drawn as continuous curves, with the
value of $n$ running from 3 to 19 (scanning left to right the curves at the top of the figure).}
\label{fig-errb1}
\end{figure}

Because of symmetry, diagonal Jacobi entries are constant and equal within numerical precision for both $\mu$ and $\mu_n$, permitting us to focus on outdiagonal ones, denoted by $b_l$ and $b^n_l$, respectively. In Figure \ref{fig-errb1} we plot $\Delta^b_{n,j}$, the average difference between the first $j+1$
entries of $J(\mu)$ and $J(\mu_n)$:
\beq
  \Delta^b_{n,j} := \frac{1}{j+1} \sum_{l=0}^j |b_l-b^n_l|.
\nuq{eq-neps0}
In the figure, drawn in doubly logarithmic scale, we observe a region where curves are approximately linear and parallel. In this range, at fixed $j$, the difference between the Jacobi matrix entry of $\mu$ and that of $\mu_n$ vanishes exponentially fast in $n$, until it attains machine precision. This region extends unlimitedly (within the reach of these experiments) to large values of $j$: this experimental manifestation of  stability proves that the procedure of the previous section defines an algorithm to compute the Jacobi matrices of balanced IFS measures.

Let us now examine quantitatively the performance of this algorithm, by considering the absolute error in the computation of $J(\mu)$.
Let us fix a threshold $\epsilon$ and let us find the rank of the largest truncation of $J(\mu)$, that is approximated within $\epsilon$ by $J(\mu_n)$. We call this rank $N_\epsilon(n)$. It can be formally defined via a quantity that will be repeatedly used in the following: the {\em $\epsilon$--coincidence range}, $\Lambda_\epsilon(J,J')$, of two Jacobi matrices $J$ (with entries $a_j,b_j)$ and $J'$ (with entries $a'_j,b'_j)$ is the integer number
\beq
  \Lambda_\epsilon(J,J') := \max \{l \mbox{  s.t. } |a_j-a'_j| \leq \epsilon, \; |b_j-b'_j| \leq \epsilon, \; 0 \leq j \leq l \}.
  \label{eq-neps3}
  \end{equation}
Using this quantity, we can write that $N_\epsilon(n) = \Lambda_\epsilon(J(\mu),J(\mu_n))$.
Letting the threshold take the value $\epsilon= 10^{-8}$  yields the values reported in Table \ref{tab-1}. Since the number of atoms in $\mu_n$ is $N=2^n$, these data can be well described by a power--law of the kind $N_\epsilon \sim A(\epsilon) N^\beta$, where $\beta$ can be estimated as $\beta \simeq .912$.
Let us now go back to Figure \ref{fig-errb1}, where $N_\epsilon(n)$ is approximately the abscissa of the intersection of the plotted curves with an horizontal line at ordinate $\epsilon$. In the linear range in  Fig. \ref{fig-errb1}, one also finds that $A(\epsilon) \sim C \epsilon^{\eta}$, with $\eta\simeq .286$ and $C$ a quantity independent of $\epsilon$ and $n$.
Taking into account that the computational complexity of the revised RKPW algorithm is approximately of $B N_\epsilon N$ arithmetical operations (with $B$ a small constant, see Appendix B), we can conclude that estimating the Jacobi matrix of $\mu$ of size $N_\epsilon$ within $\epsilon$ via this procedure has a cost that scales roughly as $\epsilon^{-\eta/\beta} N_\epsilon^{1+1/\beta} \simeq \epsilon^{-.31} N_\epsilon^{2.09}$, that is, slightly more than quadratic in the size $N_\epsilon$ and slowly increasing with respect to the precision $1/\epsilon$. %
Remark that the recursive algorithms for IFS with a
finite number of maps  \cite{cap,arxiv} require an order of $N_\epsilon^2$ operations to compute the same matrix, but this computation is exact in principle, in practice affected by a slowly increasing error \cite{arxiv}. The new algorithm \footnote{Notice that  in this section the computation of the effective size $N_\epsilon(n)$ has been effected {\em a posteriori}, using an independent knowledge of the Jacobi matrix of $\mu$. A simple technique to overcome this difficulty can be devised, mimicking Algorithm 1 of the next section. To keep the presentation contained, we do not lay down explicitly these steps, that compose Algorithm 0 in this paper.} is therefore not optimal for balanced measures, but it can be extended to equilibrium ones, for which no alternative techniques are available, except for an algebraic procedure that can be set up using a Pad\'e scheme \cite{alphonse2}.

\begin{table}
\centering
\begin{tabular}{|c|r|r|r|r|r|r|r|r|}
  \hline
  $n$ & 12 & 13 &14 &15 &16 &17 &18 &19 \\
  \hline
  $N_\epsilon(n)$  &   53 &   101 &   188 &  367 &
   693 &   1281 &   2374 &   4403 \\
  \hline
\end{tabular}
\caption{Rank $N_\epsilon(n)$ versus $n$ at $\epsilon = 10^{-8}$, for the IFS of Example \ref{examp-cantor}, Experiment 1. \label{tab-1}}
\end{table}

\section{The Jacobi matrix of the equilibrium measure on a set of intervals}
\label{sec-jeq}

The techniques of the previous two sections can also be applied to the double--index sequence of measures $\nu_n^G$.
When the order $n$ is kept fixed and $G$ tends to infinity, $\nu_n^G$ tends weakly to $\nu_n$, the equilibrium measure on $E_n$, and its Jacobi matrix tends to $J(\nu_n)$: numerically, this problem is interesting in itself, and we want to investigate the nature of this convergence. We set up the following algorithm:
\begin{itemize}
\item[] {{\bf Algorithm 1. Computing the Jacobi matrix $J(\nu_n)$.} \\
 {\bf Input}: The number of Gaussian points $G$, the threshold $\epsilon$, the IFS parameters, the order $n$ of the IFS construction and a small positive integer $g$.
 \\ {\bf Output}: the (truncated) Jacobi matrix $J(\nu_n)$, of rank $N_\epsilon(n,G)$ and absolute component-wise error less than $\epsilon$.}
\item[1:] Compute points and weights for $\nu_n^G$, eq. (\ref{eq-meas20}), solving the non-linear system of eqs. (\ref{meas5}), as described in ref. \cite{dolo} and in Appendix A.
\item[2:] Compute $J(\nu_n^G)$ via the revised RKPW algorithm, Appendix B.
\item[3:] Decrease $G$ by the amount $g$ and redo steps [1] and [2].
\item[4:] Compute $N_\epsilon(n,G) := \Lambda_\epsilon(J(\nu_n^G),J(\nu_n^{G-g}))$ as in eq. (\ref{eq-neps3}).
\end{itemize}

A few remarks are in order.
The above algorithm works for any finite set of intervals,
and not only for IFS intervals.
The cardinality of the atoms of $\nu^n_G$ is $G \times M^n$: a set of $G$ Gaussian points is used for each interval $E_n^j$, $j=1,\ldots,M^n$. This number is typically much larger than $N_\epsilon(n,G)$, so that the improvement of the RKPW algorithm presented in Appendix B becomes crucial.
Finally, the quantity $N_\epsilon(n,G)$  is a numerical estimate of the true {approximation range} $\Lambda_\epsilon(J(\nu_n^G),J(\nu_n))$.
Using this quantity we will say that the Jacobi matrix $J(\nu_n)$ of size $N_\epsilon(n,G)$ has been {\em effectively computed within $\epsilon$}.
We start now by testing the performance of this algorithm in the Julia set example described above.

{\em Numerical experiment 2. Consider the IFS in Example \ref{exa-julia} with $\lambda = 2.1$ and apply algorithm 1 with $g=1$.}

\begin{figure}
\centerline{\includegraphics[width=.6\textwidth, angle = -90]{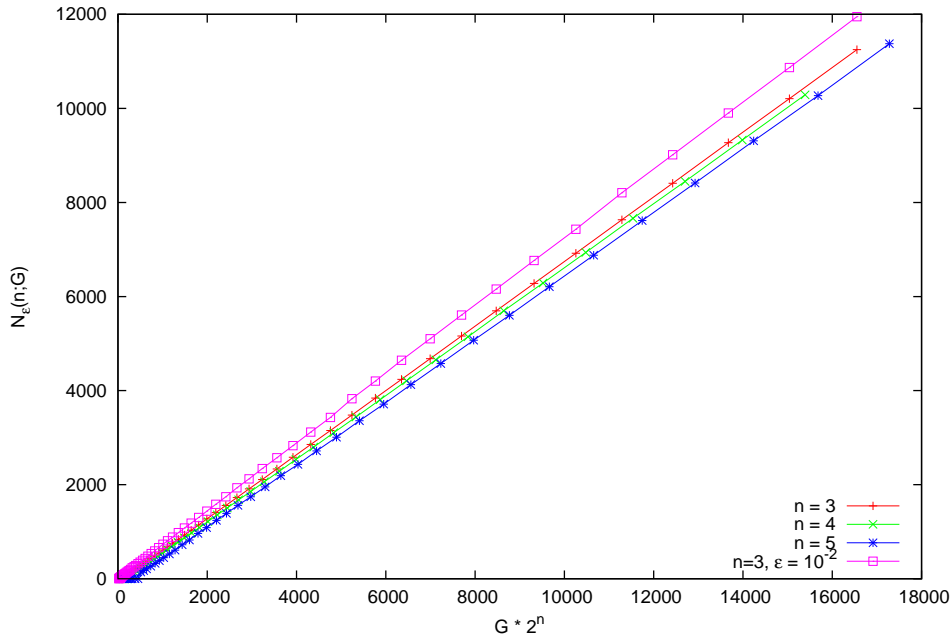}}
\caption{Size of the effectively computed Jacobi matrix, $N_\epsilon(n,G)$, versus $G \times 2^n$, for the IFS in Example \ref{exa-julia}, Experiment 2. The first three sets of data are computed with $\epsilon=10^{-12}$, the last is the case $n=3$, $\epsilon=10^{-2}$.}
\label{fig-julia1}
\end{figure}

In Figure \ref{fig-julia1} we plot the rank of the effectively computed Jacobi matrix, $N_\epsilon(n,G)$, as a function of the total number of Gaussian points, $G \times 2^n$. Data for $n=3,4,5$ and $\epsilon=10^{-12}$ show a similar linear behavior, $N_\epsilon(n,G) \sim A \; G \times 2^n$. The slope $A$ of these linear laws varies little from one case to the next, with a slightly decreasing trend. For instance, data for $n=3$ can be fitted by $A \simeq .675$, while the slope for the case $n=5$ is $A \simeq .647$.
In addition, in Figure \ref{fig-julia1} we have also plotted data for $n=3$, $\epsilon=10^{-2}$: they lie approximately on a line of slope $A \simeq .725$.

Consider now a horizontal line in the plot, going from left to right at fixed ordinate, say $N_\epsilon = H$. This line encounters the data for $n=3$, $\epsilon=10^{-2}$, at a certain value of the abscissa. This is the number of Gaussian points required to compute the truncated Jacobi matrix $J(\nu_n)$, of rank $H$, within maximum component-wise error $\epsilon$. Still moving to the right, one rapidly encounters the data for $n=3$, $\epsilon=10^{-12}$, at a new value of the abscissa. In the segment between these two points of intersection the error decreases by orders of magnitude---quickly reaching the maximum precision attainable.
Indeed, in Algorithm 1 we cannot set a threshold $\epsilon$ smaller than what allowed by the algorithm and the machine on which it is implemented. The point at which the algorithm hits the maximum precision can be found by the break-up of the increase of the approximation range, $N_\epsilon(n,G)$.

{\em Numerical experiment 3. Compute, in the same case of Experiment 2, the number of Gaussian points $\tilde{G}$ at which a further increase of $G$ does not lead to an increase of the approximation range $N_\epsilon(n,G)$.}

\begin{figure}
\centerline{\includegraphics[width=.6\textwidth, angle = -90]{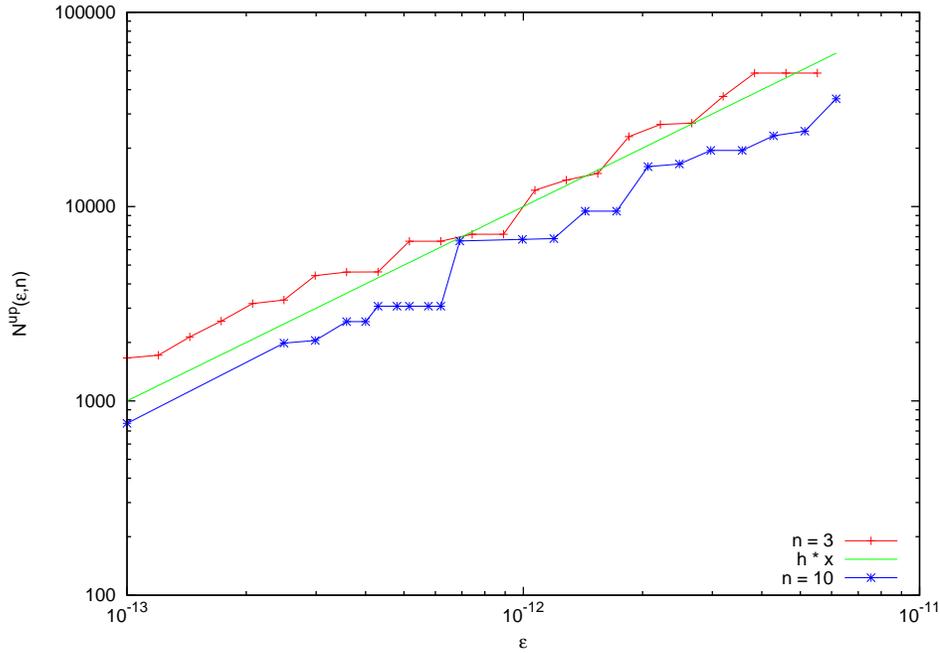}}
\caption{Values of $N^{up}(\epsilon,n)$, the maximum size of the Jacobi matrix $J(\nu_n)$ that can be computed with absolute error less than $\epsilon$, versus $\epsilon$, with $n=3$ and $n=10$. The continuous line has equation $N^{up}(\epsilon,n) = h \epsilon$, $h=10^{-16}$. Experiment 3, IFS of Example \ref{exa-julia}.}
\label{fig-julia2}
\end{figure}

Using the value $\tilde{G}$ we can find the maximum rank reachable by the algorithm, $N^{up}(\epsilon,n)=N_\epsilon(n,\tilde{G})$. In Figure \ref{fig-julia2} we plot this quantity versus $\epsilon$, for $n=3$ and $n=10$. No effort was made to have particularly clean data, since our goal here is just to estimate the law of error growth. Both sets of data are consistent with a linear dependence of $N^{up}(\epsilon,n)$ with $\epsilon$. This fact can be turned around, to imply that the minimal numerical error in the determination of the $j$-th line of the Jacobi matrix $\{(a_j$, $b_{j})\}_{j \in \bf N}$ grows linearly with $j$, when running Algorithm 1. We deem this to be an optimal result.

\section{The Jacobi matrix of the equilibrium measure on a IFS attractor}
\label{sec-jeqattra}

We are now equipped with the building blocks to compute the Jacobi matrix $J(\nu_E$), the main goal of this work. We will combine the theory of Sect. \ref{sec-equil} with the numerical techniques of Sect. \ref{sec-gragg} and \ref{sec-jeq}. Our approach is to increase alternatively the two indices of $\nu_n^G$. This is implemented in the following algorithm:

\begin{itemize}
\item[]
{{\bf Algorithm 2. Computing the Jacobi matrix $J(\nu_E)$.} \\
 {\bf Input}: The IFS parameters, the maximum order $\bar{n}$, the initial number of Gaussian points $G_0$, the threshold $\epsilon$, an increase ratio $\eta>1$.
 \\ {\bf Output}: the (truncated) Jacobi matrix $J(\nu_E)$, of rank $H(\epsilon,\nu_E)$ and absolute component-wise error less than $\epsilon$.}
 \item[1:] Let $n=2$,  $G=G_0$
\item[2:] Compute the truncated Jacobi matrix $J(\nu_{n-1}^{MG})$ and its effective size at approximation $\epsilon$, $N_\epsilon(n-1,M G)$, via Algorithm 1.
\item[3:] Similarly, compute the truncated Jacobi matrix $J(\nu_n^G)$, and $N_\epsilon(n,G)$. Let $N = \min (N_\epsilon(n-1,M G),N_\epsilon(n,G))$.
\item[4:] Compute  $H_\epsilon(n) = \Lambda_\epsilon(J(\nu_{n-1}^{MG}),J(\nu_n^G))$.
\item[6:] If $H_\epsilon(n)$ is equal to $N$, increase $G$ to $\eta G$ (or to $G+1$ when this is larger) and loop back to step 2. Otherwise continue.
\item[7:] If $n$ is smaller than $\bar{n}$ increase $n$ by one and loop back to step 2. Otherwise stop, since the truncated Jacobi matrix $J(\nu_{\bar{n}}^G)$ of size $H(\epsilon,\nu_E) = H_\epsilon(\bar{n})$ is an $\epsilon$--approximation of $J(\nu_E)$.
\end{itemize}

For any $n$, starting from the lowest, this algorithm computes finite truncations of $J(\nu_n)$ and $J(\nu_{n-1})$, accurate within a threshold $\epsilon$, with $G$ and $MG$ Gaussian points, respectively (recall that $M$ is the number of IFS maps). While comparing these matrices component--wise, it raises $G$ in order to find the maximum rank in which they coincide, again within the threshold $\epsilon$. This rank naturally increases with the order $n$.
We now investigate the performance of this algorithm in the two examples described above.

\subsection{Julia set equilibrium measure}

Recall that we have at our disposal simple recursion relations \cite{belli}, reproduced in Sect. \ref{sec-ifs}, for the Jacobi matrix associated with
Julia set measures, that can serve to gauge the precision of Algorithm 2.

{\em Numerical experiment 4. Run Algorithm 2 on Example
\ref{exa-julia}, with $\lambda=2.1$ and $\epsilon=10^{-10}$.}

Table \ref{tab-4} is a report of the quantities characterizing this experiment. Observe first that, when using the same number of {\em total} Gaussian points, $J(\nu_n)$ can be effectively computed 
of roughly twice the rank of $J(\nu_{n-1})$. The entries of these matrices differ less than $\epsilon$ up to the index denoted $H_\epsilon(n)$: this number is much smaller than both the preceding ranks, $N_\epsilon(n-1,M G)$ and $N_\epsilon(n,G)$, but like these it increases geometrically with $n$.
The last column is the ``true'' effective size
$Y_\epsilon(n) := \Lambda_\epsilon(J(\nu_n^G),J(\nu_E))$, in which $J(\nu_E)$ has been computed ``exactly'' by the recursion relations in Section \ref{sec-ifs}. It is remarkable that the value of $Y_\epsilon(n)$ is equal to $H_\epsilon(n+1)$, and that this value is $2^n -1$.

\begin{table}
\centering
\begin{tabular}{|r|r|r|r|r|r|}
  \hline
  $n$ & $G$ & $N_\epsilon(n-1,M G)$ & $N_\epsilon(n,G)$ & $H_\epsilon(n)$ & $Y_\epsilon(n)$  \\
  \hline
    2    &      13    &      11      &     7    &       1    &       3 \\
    3    &      13    &       7    &     24     &     3     &     7 \\
     4   &      13    &     24     &    55     &    7      &  15\\
      5  &       13   &      55    &    112    &     15    &     31\\
      6   &       13   &      112   &      224   &       31  &        63\\
     7    &     13     &   224   &   451    &  63     &    127\\
   8    &   13     &   451     & 943     &  127     & 255\\
  9    &    12     &   511    &  1567    &    255    &     511\\
  10   &      12   &    1567   &     3135   &     511   &     1023\\
\hline
\end{tabular}
\caption{Numerical results of Experiment 4, for Example \ref{exa-julia}. See text for details. \label{tab-4}}
\end{table}

Some of the phenomena observed in the previous experiment are clearly typical of Julia sets, while others are more general. Observe that $H_\epsilon(n)$, being a coincidence of the matrices $J(\nu_{n-1})$ and $J(\nu_{n})$ is to be regarded as an estimate of the number of the ``correct'' components  of the former, rather than the latter, a fact confirmed by the last column of Table \ref{tab-4}. It is then possible to obtain a better estimate of the number of effectively computed entries of $J(\nu_E)$ at level $n$ by the extrapolation rule $\hat{H}_\epsilon(n) = H^2_\epsilon(n)/H_\epsilon(n-1)$.

\begin{figure}
\centerline{\includegraphics[width=.6\textwidth, angle = -90]{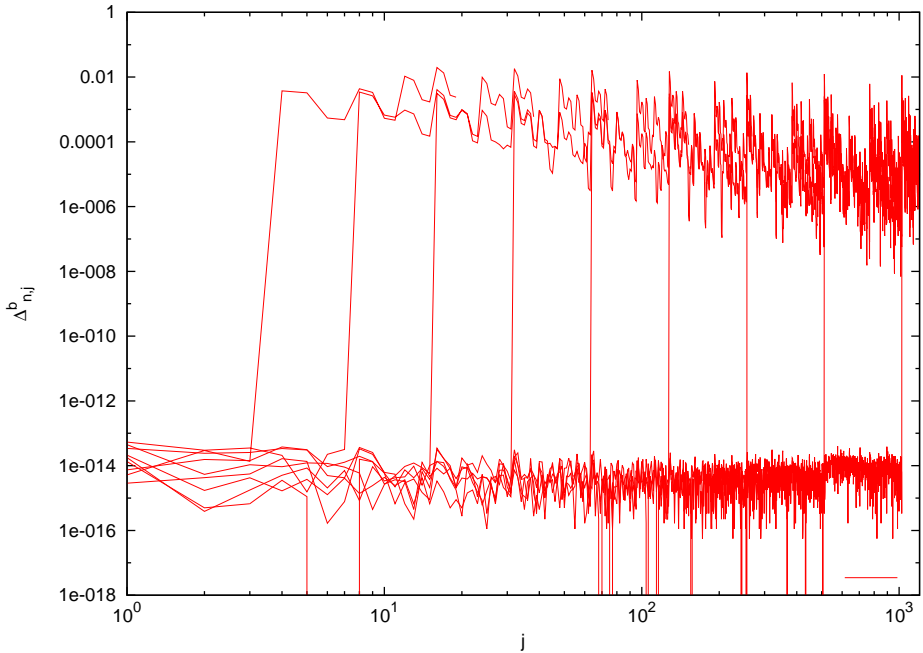}}
\caption{Absolute errors
$\Delta^b_{n,j} := |b_j(\nu_n^G) - b_j(\nu_E)|$ versus $j$, in the computation of the Jacobi matrix of Example \ref{exa-julia}, Experiment 4, with $\epsilon=4 \times 10^{-13}$. Data for $n=2$ to $n=10$ can be recognized from their sharp rise at geometrically increasing values, $j=Y_\epsilon(n)$.}
\label{fig-julia3}
\end{figure}

Experiment 4 can also be used to gauge the precision of the algorithm, and the error growth. To do this, we have
re-run it with $\epsilon=4 \times 10^{-13}$, a much smaller value than the threshold used in Table \ref{tab-4}, and we have computed the absolute error $\Delta^b_{n,j} := |b_j(\nu_n^G) - b_j(\nu_E)|$, for a range of Jacobi entries of index $j$ that exceeds $Y_\epsilon(n)$. Figure \ref{fig-julia3} shows that
$\Delta^b_{n,j}$ is roughly constant in the full range $j \leq Y_\epsilon(n)$ and then skyrockets when this value is surpassed. This also implies that $Y_\epsilon(n)$, at fixed $n$, is approximately constant in a large range of $\epsilon$ values.

\subsection{Ternary Cantor set equilibrium measure}

This example is numerically more demanding than the Julia set case we have just examined: a larger number of Gaussian points are required to compute the Jacobi matrix of $\nu_n$ and, on top of that, convergence to the limit Jacobi matrix, when increasing $n$, is slower and does not follow the simple pattern observed in Figure \ref{fig-julia3}.
Therefore, we adopt a simplified version of Algorithm 2. We start by fixing a threshold $\epsilon$ and a size $N$. We then run Algorithm 1 at increasing numbers of Gaussian points $G$, until the effective rank of $J(\nu_n)$, $N_\epsilon(n,G)$, is larger than $N$. We do this for a range of $n$ values. Finally, we compare the resulting Jacobi matrices $J(\nu_n)$, looking for convergence of their entries. This procedure is used in the following experiment.

{\em Numerical experiment 5. In the case of Example \ref{examp-cantor}, let $N=65540$, $\epsilon = 10^{-8}$, $\bar{n}=15$. Run algorithm 1 raising $G$ until  $N_\epsilon(n,G) \geq N$, for $n=1,\ldots,\bar{n}$.}

\begin{figure}
\centerline{\includegraphics[width=.6\textwidth, angle = -90]{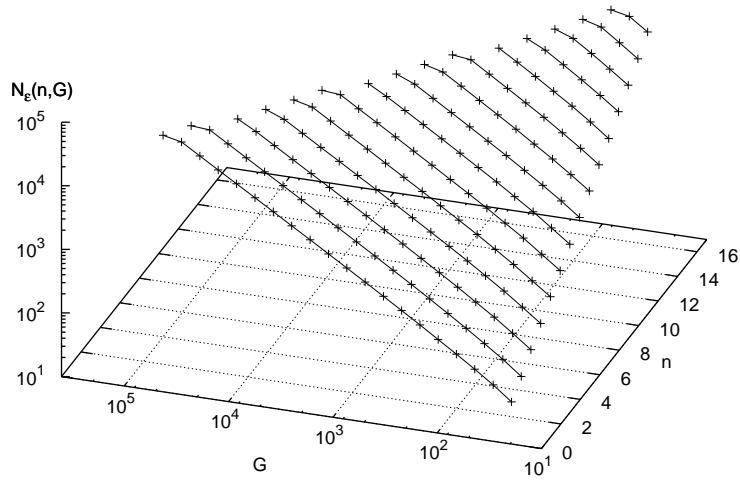}}
\caption{Values of $N_\epsilon(n,G)$ versus $n$ and $G$, for Example \ref{examp-cantor}, Experiment 5.}
\label{fig-cantor1b}
\end{figure}

\begin{figure}
\centerline{\includegraphics[width=.6\textwidth, angle = -90]{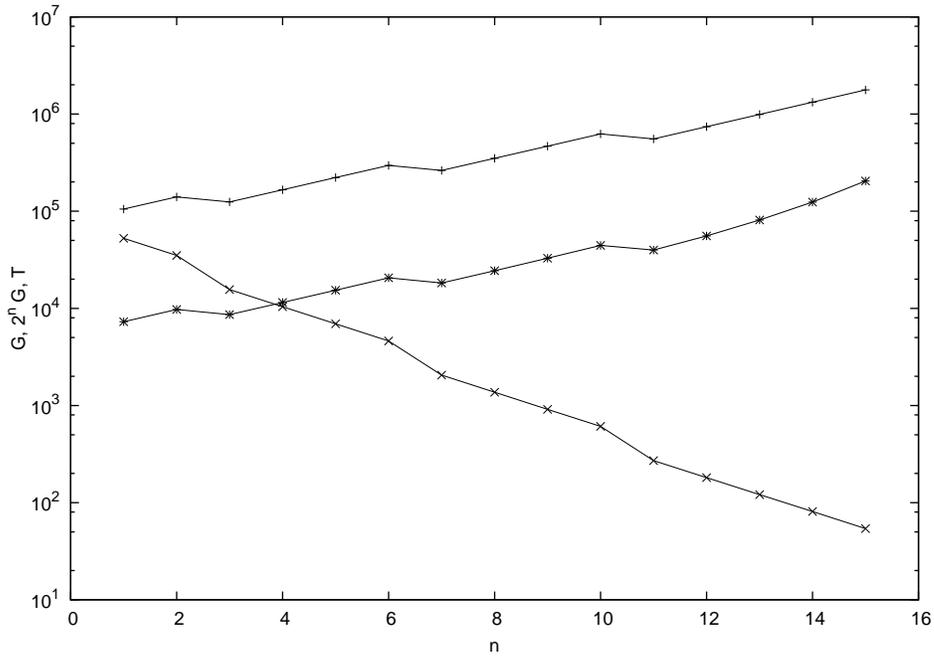}}
\caption{Number of Gaussian points $\tilde{G}$ required to compute $J(\nu_n)$ so that
$N_\epsilon(n,\tilde{G}) \geq 65,540$, with $\epsilon = 10^{-8}$ (decreasing data, crosses); total number of  Gaussian points, $2^n \tilde{G}$ (top curve, pluses) and
total execution time $T$ (seconds) (stars). Example \ref{examp-cantor}, Experiment 5.}
\label{fig-cantor2}
\end{figure}

The sequence of values of $N_\epsilon(n,G)$ computed by the algorithm are reported in Figure \ref{fig-cantor1b} versus $n$ and $G$. At fixed $n$, this figure displays convergence of the Jacobi matrices $J(\nu_n^G)$ to $J(\nu_n)$: the rank of the $\epsilon$-approximation grows linearly, as in Fig. \ref{fig-julia1}, when increasing the number of Gaussian points. When varying $n$ and $G$ at the same time, increasingly larger truncations of $J(\nu_E)$ can be computed.
The algorithm stops when $N_\epsilon(n,G)$ is larger than $N$. The number of Gaussian points at which this is achieved, $\tilde{G}$, is tabulated as a function of $n$ in Table \ref{tab-3} and is plotted in Figure \ref{fig-cantor2}. We observe that $\tilde{G}$ decreases exponentially with $n$. Yet, recalling that each interval $E_n^i$ requires $\tilde{G}$ Gaussian points, the {\em total} number of Gaussian points, $2^n \tilde{G}$, also plotted in Figure \ref{fig-cantor2}, increases at a (moderate) exponential rate, roughly equal to $(1.25)^n$. This increase is a manifestation of the difficulty of the problem, and impacts its computational complexity, seen in the third curve in the figure, where it is measured as the total number of cpu seconds required to compute all data points in a string at fixed $n$ in Fig. \ref{fig-cantor1b}, on a 36 processors parallel machine.

\begin{table}
\centering
\begin{tabular}{|c|r|r|r|r|r|r|r|r|c|r|}
  \hline
  $n$ & 1 & 2 & 3 & 4 & 5 & 6 & 7 & 8 & \ldots & 15 \\
  \hline
  $\tilde{G}$  &
  52597  & 35065   &      15585   &       10390   &      6927
                 &     4618  &     2053  &     1369  &  \ldots & 54  \\
  \hline
\end{tabular}
\caption{Number of Gaussian points $\tilde{G}$ required to compute an $\epsilon$--approximation of $J(\nu_n)$, so that
$N_\epsilon(n,\tilde{G}) \geq 65,540$, with $\epsilon = 10^{-8}$. Example \ref{examp-cantor}, Experiment 5.\label{tab-3}}
\end{table}

\begin{figure}
\centerline{\includegraphics[width=.6\textwidth, angle = -90]{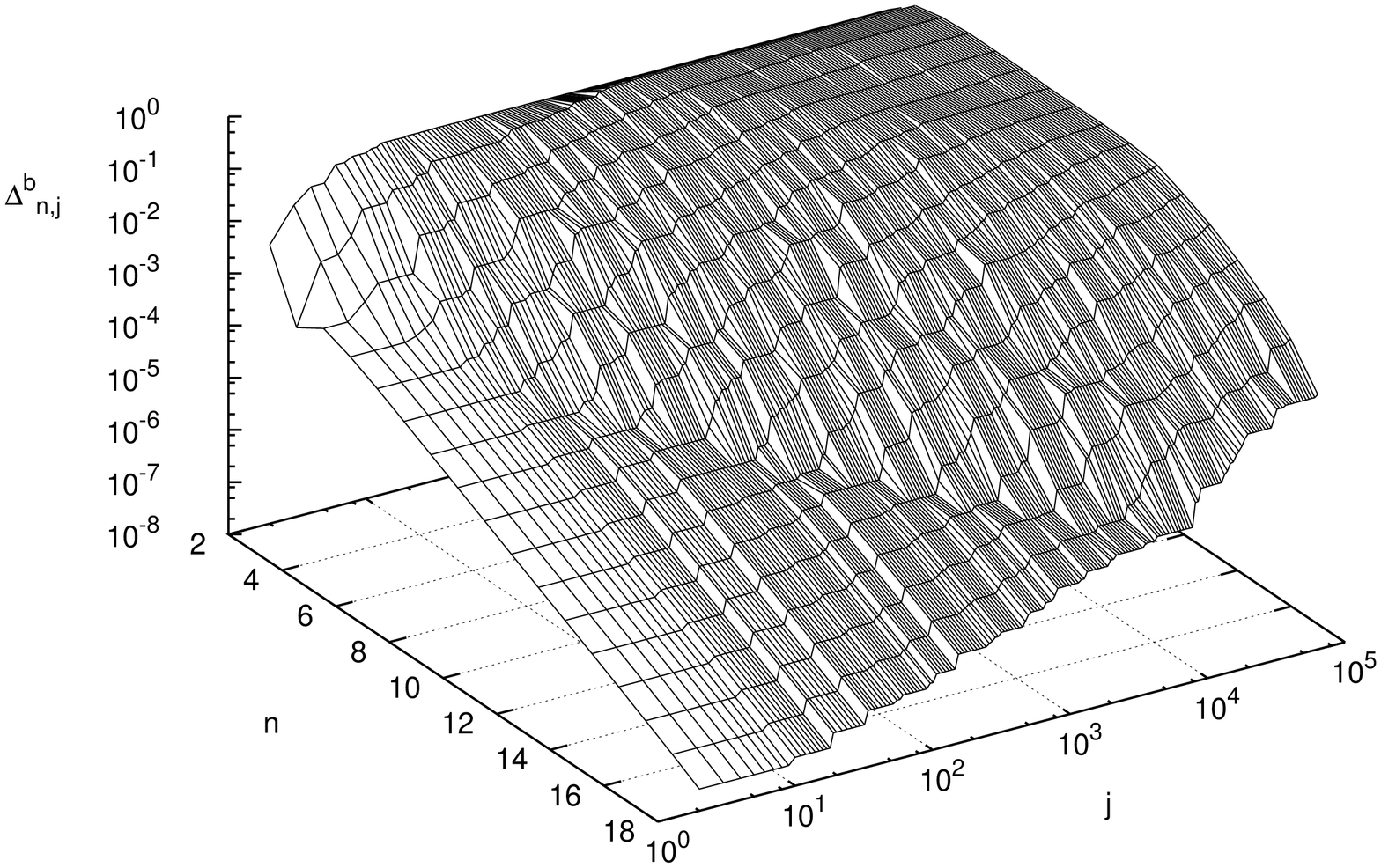}}
\caption{Absolute errors
$\Delta^b_{n,j}$, defined in the text, in the determination of the Jacobi matrix of Example \ref{examp-cantor}, Experiment 5.}
\label{fig-cantor1}
\end{figure}

\begin{figure}
\centerline{\includegraphics[width=.6\textwidth, angle = -90]{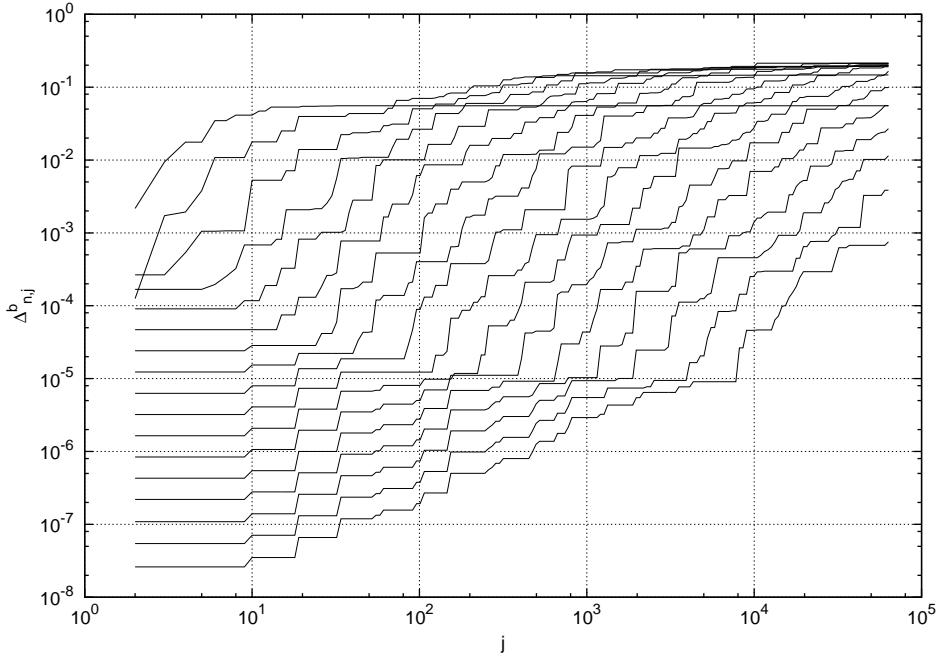}}
\caption{Absolute errors
$\Delta^b_{n,j}$, defined in the text, in the determination of the Jacobi matrix of Example \ref{examp-cantor}, Experiment 5. Curves for $n=2$ to $n=18$ are displayed,
and are monotonically ordered from top to bottom at $j=10$.}
\label{fig-cantor1bx}
\end{figure}

Finally, we define the absolute errors $\Delta^b_{n,j} := \max \{|b_l(\nu_n^{\tilde{G}}) - b_l(\nu_{n-1}^{\tilde{G}})|, \; 1 \leq l \leq j \}$.
The quantity $\Delta^b_{n,j}$, being a maximum over the first $j$ components, is an estimate of the absolute error in the determination of the entries of the matrix $J(\nu_E)$, up to rank $j$, via $J(\nu_{n-1}^{\tilde{G}})$. In Figure \ref{fig-cantor1} we plot this quantity versus $j$ and $n$. Comparison with the analogue Figure
\ref{fig-julia3} for Example \ref{exa-julia} reveals that convergence is here of a different kind: increasing $n$ at fixed $j$ yields exponential convergence, but only after a value of $n$ that increases with $j$: Figure \ref{fig-cantor1} displays in fact interesting regions with different scaling properties. A two--dimensional plot of the same data, Figure \ref{fig-cantor1bx}, in which values with the same $n$ are plotted as a line versus $j$, reveals that, for the largest plotted case, $n=18$, one can safely assume that the truncated Jacobi matrix of $J(\nu_E)$, of 65,400 entries, is estimated by $J_{18}^{\tilde{G}}$ with an error smaller than one part in a thousand.

It is then apparent that the Cantor set in Example \ref{examp-cantor}  is a more demanding test than the Julia set in Example \ref{exa-julia}. Nonetheless, the technique described in this section is capable of computing large truncations of its Jacobi matrix. As a matter of facts, the most sensitive and time--consuming step in this computation seems to be the determination of the roots $\{\zeta_i\}$ of the set of equations (\ref{meas5}).

\section{Root asymptotics of orthogonal polynomials}
\label{sec-capa}

It is well known that balanced and equilibrium measures on the ternary Cantor set $E$,  as well as the finite sets of intervals $E_n$, are regular, in the sense described in \cite{stahl}. Let $\sigma$ be any of these measures, and let $\nu$ be the measure $\nu_E$ or $\nu_n$, accordingly.
Since these sets have positive capacity, empty interior and connected complement, regularity means that the limit (\ref{pote0}) exists, and is equal to the equilibrium measure on $E$ ({\em regularity of the zero distribution}) and that, equivalently, the limit of $\log |p_j(\sigma;z)|^{1/j}$ is the Green's function ${g(E;z)}$, (or  ${g(E_n;z)}$, in the finite $n$ case), defined in eq. (\ref{pote8}), locally uniformly outside the convex hull of $E$ ({\em regular $j$-th root asymptotic behavior}): see thms 3.1.1 and 3.1.4 in \cite{stahl}. The last relation involves the modulus of the orthogonal polynomials.
It is easy to extend the j-th root behavior to the complex {\em monic}\footnote{ {\em i.e} those for which the coefficient of $z^j$ is one, used here for simplicity. Recall that these polynomials satisfy the recursion relation $P_{j+1}(\sigma;z) = (z-a_j) P_j(\sigma;z) - b_j^2 P_{j-1}(\sigma;z)$.}
orthogonal polynomials $P_j(\sigma;z)$.

Since $P_j(\sigma;z) = \prod_{l=1,\ldots,j} (z-\xi^j_l)$, regularity of the zero distribution implies that the following limit exists:
 \beq
   \lim_{j \rightarrow \infty} \frac{1}{j} \log(P_j(\sigma;z)) =  \int \log (z-s)   \; d \nu(s) :=  - \Psi(\nu;z),
 \label{eq-capa1}
 \eeq
with the principal determination of the logarithm, when $z$ does not belong to the 
support of $\sigma$. The integral in the above equation defines the complex potential $\Psi(\nu;z)$.  Its real part is the electrostatic potential, $V(\nu;z)$, eq. (\ref{pote1}), so that the real part of eq. (\ref{eq-capa1}) conveys the conventional meaning of regular root asymptotics. In the finite $n$ case, the function $\Psi(\nu_n;z)$, is the same as the {\em complex} Green's function of Widom (\cite{widom} eq. 14.1), modulo the constant  $\log(\mbox{Cap}(E_n))$:
 \beq
       \Psi(\nu_n;z) - \log(\mbox{Cap}(E_n)) = \int_{\alpha_1}^z Z(\zeta;s)/\sqrt{Y(s)} \; ds.
 \label{eq-capa1c}
 \eeq
Introducing the polynomial ratios $\rho_j(\sigma;z) := P_j(\sigma;z)/P_{j-1}(\sigma;z)$ in eq. (\ref{eq-capa1}) and taking real and imaginary parts leads to:
 \beq
  \lim_{j \rightarrow \infty}\frac{1}{j} \sum_{i=1}^j \log |\rho_i(\sigma;z)|=-\Re \Psi(\nu;z),
 \label{eq-capa2}
 \eeq
and
  \beq
  \lim_{j \rightarrow \infty}\frac{1}{j} \sum_{i=1}^j \arg (\rho_i(\sigma;z)) =-\Im\Psi(\nu;z).
 \label{eq-capa3}
 \eeq

This procedure reveals the nature of a Lyapunov exponent \cite{johnson,russell} for the real part of the complex potential, eq. (\ref{eq-capa2}); it also makes clear that weaker requirements are needed for root asymptotics, with respect to ratio asymptotics---this latter being the existence of the limit of $\rho_j$ \cite{stahl,lubin}.

Let us now come to the imaginary part, eq. (\ref{eq-capa3}).
The arguments in this equation are assumed to lie in the interval $(0,2 \pi)$.
When $z = x+iy$ tends to the real axis from the above, the left hand side of this equation becomes the conventional rotation number of the theory of discrete Sturm--Liouville operators \cite{souil,johnson}.
Coherently, in this limit the right hand side $\Im \Psi(\nu;z)$ tends to $-\pi \nu([x,\infty))$. Therefore, the imaginary part of the complex potential $\Psi(\nu;z)$ extends the rotation number to the complex plane \cite{johnson} and the left hand side of
eq. (\ref{eq-capa3}) gives a practical means to compute it.

We want now to investigate the numerical implications of eqs. (\ref{eq-capa2},\ref{eq-capa3}).
The function $\Psi(\nu;z)$ is an integral with respect to $\nu$, eq. (\ref{eq-capa1}), that can be numerically estimated by a limit procedure, quite analogous to that of Sects. \ref{sec-jeq},\ref{sec-jeqattra}, of the corresponding integrals with the discrete measures $\nu_n^G$. Moreover, using the Jacobi matrices of $\nu_n$ and $\nu_E$, we can also compute the left hand side of eqs. (\ref{eq-capa2},\ref{eq-capa3}) and verify experimentally the rate at which convergence takes place in root asymptotics.

{\em Numerical experiment 6. In the case of Example \ref{examp-cantor}, Experiment 5, choose $\sigma=\nu_n$ and compute the sequence of real parts $g_{n,j}(z) := \log(|P_j(\nu_n;z)|)/{j}$ using the recursion relation $\rho_{i+1} = (z-a_i)\rho_i - b_i^2 / \rho_{i-1}$ for the left hand side of eq. (\ref{eq-capa2}). Using the discrete measures $\nu_n^G$ in the integral (\ref{eq-capa1}), also compute the potential $V(\nu_n;z)$. }\\

\begin{figure}
\centerline{\includegraphics[width=.6\textwidth, angle = -90]{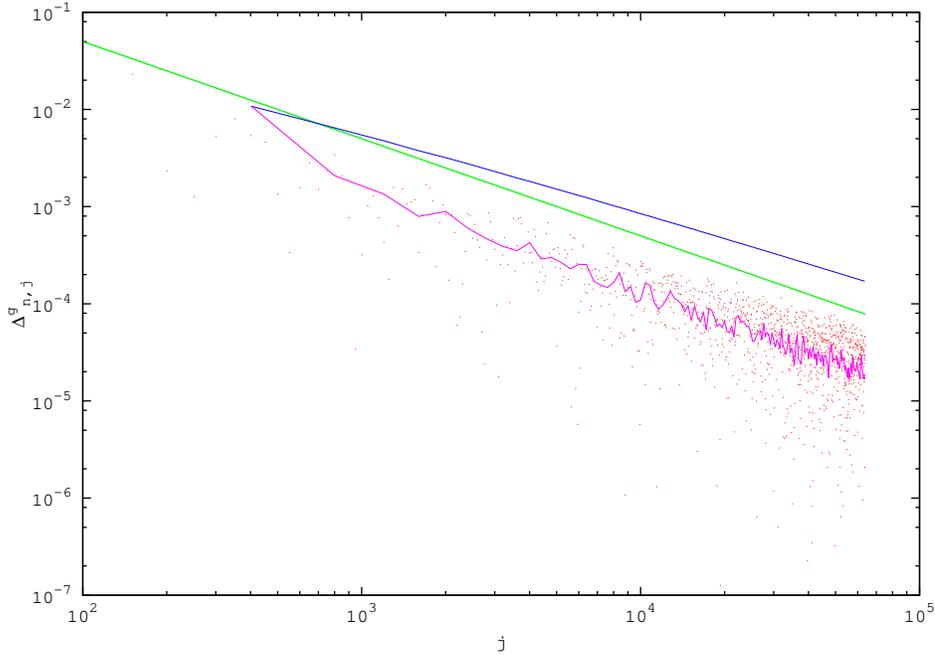}}
\caption{Absolute difference
$\Delta^g_{n,j}$ defined in the text, in the convergence of the root asymptotics for monic orthogonal polynomials $P_j(\nu_{n};z)$, plotted at equally spaced values of $j$, with $n=18$, $z=\frac{1}{4} + 5 \cdot 10^{-6}$. Example \ref{examp-cantor}, Experiment 6. We plot: the raw values $g_{n,j}$ (red dots), the moving C\'esaro averages of $g_{n,j}$ over the interval $[j,j+L]$, $L=400$ (magenta line) and the C\'esaro average over the interval $[1,j]$ (blue line). Also drawn is the function $h(x)= 5/x$ (green line). }
\label{fig-capa1}
\end{figure}

Distance from the asymptotic limit can be gauged by the differences
 \beq
  \Delta^g_{n,j} = |g_{n,j}(z) -  V(\nu_n;z)|.
 \label{eq-capa2x}
 \eeq
In Figure \ref{fig-capa1} these differences are plotted versus $j$, for $n=18$, at a specific value of $z$, very close to $E_n$. We observe wide fluctuations accompanying the general decay assured by the theory.
To smooth these oscillations, we first average $g_{n,l}(z)$ over a set of $L$ integer indices, ranging from $l=j$ to $l=j+L-1$, and then compute the new difference from the exact potential, using again eq. (\ref{eq-capa2x}). Both these sets of data are consistent with power-law convergence, with exponent minus one, a fact on which we will comment momentarily.
Finally, we also perform the full C\'esaro average of $g_{n,l}(z)$ for $l$ from one to $j$. In this case we pay the smoother data by a slower convergence rate (algebraic decay with exponent $\eta \simeq -0.85$). Observe that in this last case we are indeed performing a full double C\'esaro average of the real part of the logarithm of $P_j(\nu_n;z)$. It might be that this last set of data is better suited to extrapolate $V(\nu_n;z)$ via any of the usual techniques \cite{claudemichela}. Yet, in the following we will simply use the C\'esaro average of  $L$ values $g_{n,j}(z)$, at the largest available indices, as a numerical approximation for the limit $j \rightarrow \infty$ in eq. (\ref{eq-capa2}).

\begin{figure}
\centerline{\includegraphics[width=.6\textwidth, angle = -90]{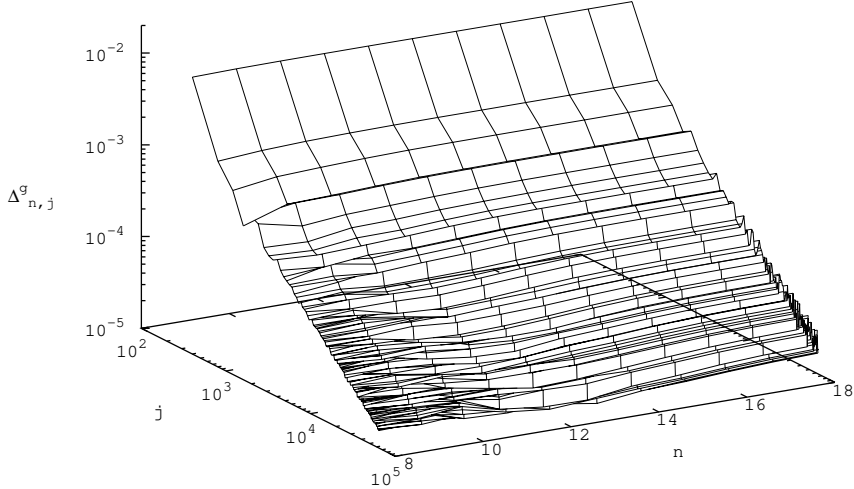}}
\caption{Absolute difference
$\Delta^g_{n,j}$ defined in the text, for the moving C\'esaro averages of $g_{n,j}$ over the interval $[j,j+L]$, $L=400$, plotted versus $n$ and $j$, at $z=\frac{1}{4} + 5 \cdot 10^{-6}$. Example \ref{examp-cantor}, Experiment 6.
 }
\label{fig-capa2}
\end{figure}

The data just shown are related to the case of $\nu_n$, $n=18$, for which the difference between $J(\nu_n)$ and $J(\nu_E)$ has been estimated in the previous section. In Figure \ref{fig-capa2} we plot the function $\Delta^g_{n,j}$, for local C\'esaro averages, as a function of both the polynomial index $j$ and the order $n$ of the hierarchical construction of the Cantor set. We observe a remarkable similarity of the curves for different $n$, which decay as $j^{-1}$, that might indicate, in the infinite $n$ limit,
an asymptotic formula of the kind $|P_j(\nu_E;z)| = B_j(\nu_E;z) \; e^{-j V(\nu_E;z)}$ with $B_j(\nu_E;z)$ a bounded function of $j$.
To prove this conjecture rigorously, one might try to use the formulae, explicit albeit involved, that exist for $P_j(\nu_n;z)$ (see for instance eq. 2.25 in \cite{franzimrn}).
In any case, the observed numerical behavior guarantees that we can compute the logarithmic potential via root asymptotics: in fact, the same analysis can be carried out with similar results for the imaginary parts in eq. (\ref{eq-capa3}). The full complex potential $\Psi(\nu;z)$ is needed in the next section.

\section{Conformal mappings in constructive function theory}
\label{sec-confor}
The technique presented in the previous section also provides an efficient computational tool to investigate a conformal mapping of interest in constructive function theory \cite{andriev},\cite{andriev2},\cite{andriev3},\cite{totik}. In fact,
the function $F(\nu_E;z)$,
 \beq
   F(\nu_E;z) =  \exp \{ \Psi(\nu_E;z) - \log (\mbox{Cap} (E)) \},
 \label{eq-conf1}
 \eeq
defined in the upper half plane ${\bf H} = \{ z : \Im (z) > 0 \}$ is univalent, and enjoys significant analytical properties, that have been related to the geometry of $E$ using conformal invariants \cite{pomme,andriev}. This function maps conformally ${\bf H}$, deprived of the real set $E$, to the exterior of the unit disk, with radial spikes corresponding to the gaps in the set $E$ \cite{andriev}.

The numerical computation of conformal maps is a classical problem \cite{henrici},
that has been solved with such a variety of techniques \cite{martin},\cite{lothar},\cite{driscoll},\cite{nick},\cite{shen} that any list of references is forced to be tentative. In our problem, we face the problem of a multiply connected region of high genus---actually, for the Cantor set this genus is infinite: this renders most of the previous techniques inefficient. We have found that using the analytical solution for the finite interval case, $E_n$, and taking it to the infinite $n$ limit, is the best approach. In addition, we want to describe in this section a further observation: not only the information encoded in the Jacobi matrix of $\nu_E$ is sufficient for this task (this ought to be trivial, since $\nu_E$ is in one--to--one relation with its Jacobi matrix), but also it yields the most efficient procedure from a computational viewpoint. The reason for this is the fast convergence of the {\em complex} root asymptotics, described in the previous section.

The potential $\Psi(\nu_E;z)$ and the conformal mapping $F(\nu_E;z)$ can be computed via the integral (\ref{eq-capa1}) and Gaussian summation with $\nu^G_n$, as done in a part of Experiment 6. In so doing, the reliable computation of each functional value $F(\nu_E;z)$ requires a summation over the large number of Gaussian points $G \times M^n$. Obviously, one can optimize $n$ and $G$ versus precision, but this burden affects any point $z$: typically, many such points are required to have an illustration of the mapping $F$---or the scaling behavior of $\Psi(\nu_E;z)$ on sequences of complex points $z$, in investigations like \cite{totik,andriev3}.
To the contrary, the results of the previous section provide us with an alternative technique, of vastly inferior complexity.
In fact, we have shown that $\Psi(\nu_E;z)$ can be reliably obtained by local C\'esaro averages of the root asymptotics of orthogonal polynomials $P_j(\nu_E;z)$. Here, Gaussian summations must be performed only once, in the construction of the Jacobi matrix $J(\nu_E)$. The computation of the root asymptotics, for any point $z$, is then characterized by an extremely low complexity, that scales linearly with the maximum size of the Jacobi matrix involved. We now describe the results of this approach in our most challenging example, the ternary Cantor set, Example \ref{examp-cantor}.

{\em Numerical experiment 7. In the case of Example \ref{examp-cantor}, compute the potential $\Psi(\nu_n;z)$ via root asymptotics, eqs. (\ref{eq-capa2},\ref{eq-capa3}), for sets of complex values $z_k=x+i y$ in ${\bf H}$, lying either on horizontal lines at fixed ordinate $y$, or on vertical lines at fixed abscissa $x$. Also via root asymptotics, compute $\mbox{Cap} (E_n)$ via $V(\nu_n;z_0)$, for $z_0 \in E_n$.}

\begin{figure}
\centerline{\includegraphics[width=.6\textwidth, angle = -90]{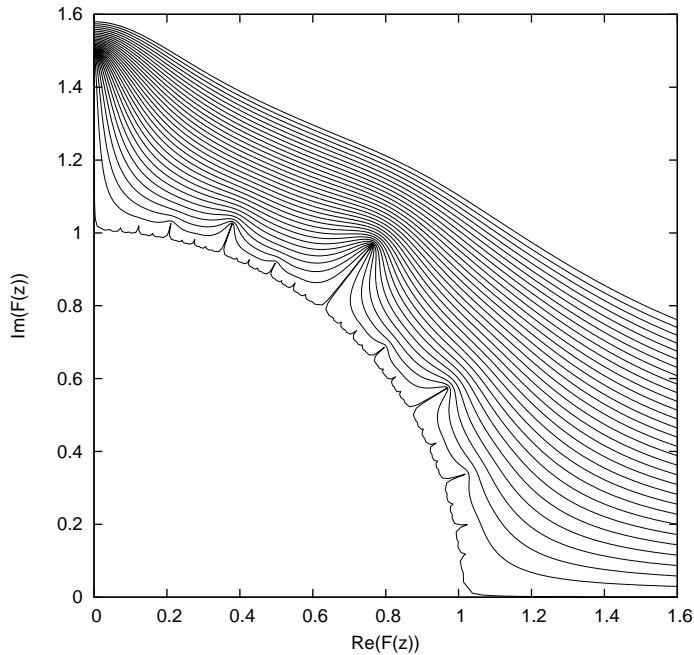}}
\caption{A Cantor Sun shining in the complex sky: Conformal mapping $F(\nu_{18};z)$ for the ternary Cantor set equilibrium measure, Example \ref{examp-cantor}, Experiment 7. The image of 31 horizontal segments $z = x + i y$, with $x$ in the interval $[.5,1.2]$ discretized in 2,000 points is drawn. The ordinates $y$ range from $5 \cdot 10^{-5}$ to $4 \cdot 10^{-2}$ over 31 intermediate values. The Jacobi matrix $J(\nu_{18})$ as been employed, with averaging in the polynomial index $j$ from 48,000 to 50,000. Within the graphical resolution of the plot, the picture is hardly distinguishable from that of $F(\nu_E;z)$.}
\label{fig-confo1}
\end{figure}

\begin{figure}
\centerline{\includegraphics[width=.6\textwidth, angle = -90]{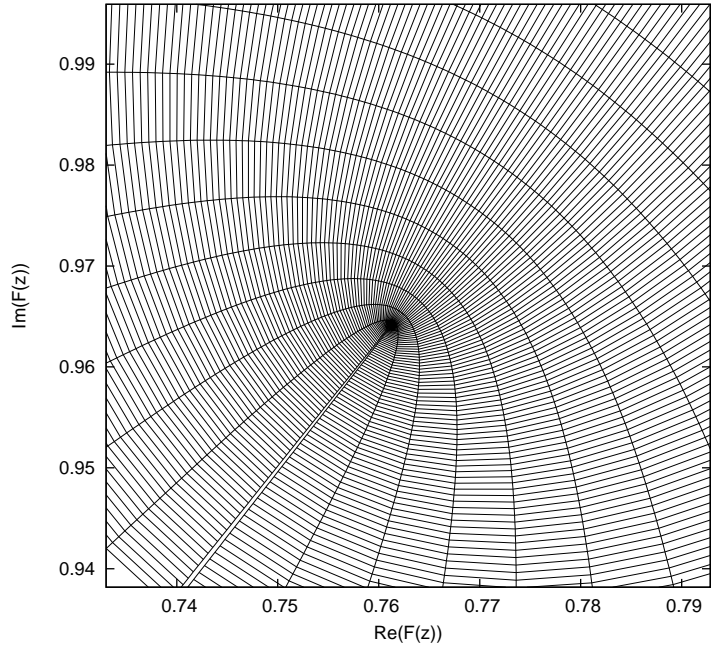}}
\caption{Conformal mapping $F(\nu_{18};z)$ for the Cantor set equilibrium measure, showing a magnification of a part of Figure \ref{fig-confo1}, also displaying the image of vertical segments.  }
\label{fig-confo2}
\end{figure}

Observe that in \cite{dolo} we have computed $\mbox{Cap} (E_n)$ via integrals of the kind (\ref{pote1}), while here we use root asymptotics also for this goal.
Plotting the values of $F(\nu;z_k)$ joined by lines gives a pictorial illustration of the properties of the conformal mapping $F$: Figure \ref{fig-confo1} displays the images of the horizontal segments, and Fig. \ref{fig-confo2} is a magnification, showing both horizontal and vertical segments, of the region nearby the tip of one of the spikes. Let us now conclude with an algorithmic interpretation of these results.

\begin{figure}
\centerline{\includegraphics[width=.6\textwidth, angle = -90]{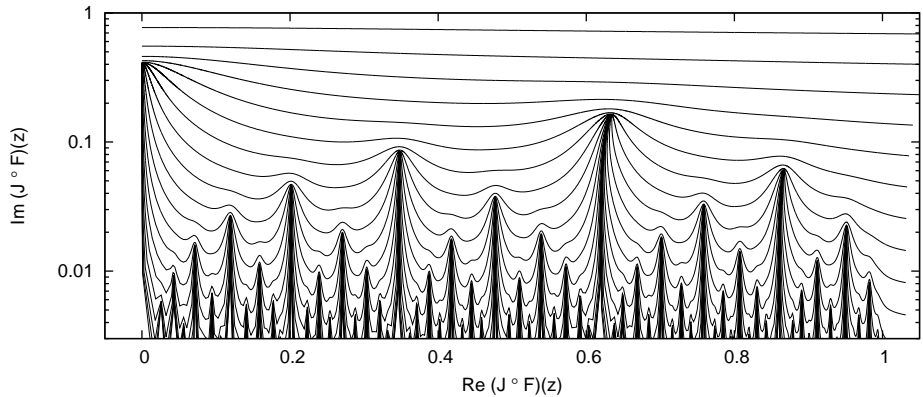}}
\caption{Conformal mapping $(J \circ F_{\nu_{18}})(z)$ for the Cantor set equilibrium measure, where the horizontal segments are taken at exponentially spaced values of $y$ ranging from $5 \cdot 10^{-6}$ to $5 \cdot 10^{-1}$.  }
\label{fig-confo3}
\end{figure}

Theoretically, it is well known that the behavior of $F(\nu_E;z)$ reveals the fine structure of the Cantor set when $x = \Re z \in E$ and $y = \Im z$ tends to zero. On the other hand, it is visually evident in Figure \ref{fig-confo1} that $F(\nu_n;z)$ reliably approximates $F(\nu_E;z)$ the more the ordinate $y$ is larger than the length of the spikes corresponding to gaps at order larger than $n$ in the Cantor set construction.  This is even more evident when composing $F(\nu_E;z)$ with the familiar Joukovsky mapping $J(z) =(z+z^{-1})/2$. The resulting conformal map $J \circ F(\nu_E;\cdot)$, from ${\bf H}$ into ${\bf H}$ minus an infinite set of slits, is conveniently seen by plotting the imaginary part in logarithmic scale, as in Fig. \ref{fig-confo3}. In this picture we can observe the hierarchical organization of slits corresponding to the countable sets of gaps of the Cantor set. This proves that the information encoded in the Jacobi matrix $J(\nu_E)$ has been correctly retrieved, and provides further evidence that the algorithm to compute such matrix (the encoding step) is experimentally stable and can reach large polynomial orders.

\section{Conclusions}

We have presented in this paper a numerical study of the equilibrium measure on attractors of Iterated Function Systems and of its orthogonal polynomials. This study employs a sequence of numerical techniques, suitably enhanced, when not explicitly designed, for IFS construction. We have described background and motivations of this research in the Introduction, so that we need only to briefly mention here a possible extension of this research, that is almost immediate. One can compute equilibrium measures and Jacobi matrices for Cantor sets of positive Lebesgue measure \cite{franzinfty,chris3}. These measures lie in between ``conventional'' absolutely continuous measures and those studied in this paper and therefore are attractive objects to analyze, theoretically and numerically.
\\

{\large \bf Acknowledgements}
I would like to thank Alphonse Magnus for many discussions and exchanges of numerical data on the equilibrium problem for Cantor sets. 
I also benefitted in the past from precious advice on the finite interval case from a sorely missed friend, Franz Peherstorfer.

\section{Appendix A: solving the nonlinear equations (\ref{meas5})}

The numerical solution of the equilibrium problem on a finite union of real intervals has been already discussed in a number of papers, and good references are \cite{nick},\cite{shen}. In these works, the polynomial $Z(z)$ in eq. (\ref{meas3}) is developed in monomials, following \cite{widom}, Sect. 14, and therefore the expansion coefficients are defined by a system of $N-1$ {\em linear} equations. Clearly, linearity is an advantage, that has lead \cite{nick} to call this the {\em trivial} direction of the Schwarz--Christoffel mapping problem. This is certainly true for problems involving few intervals.
Yet, this advantage quickly evaporates for larger numbers of intervals (of the order of the hundreds of thousands in this paper) because of the poor conditioning of the system of linear equations. A cure to this problem might be to expand $Z(z)$ in a different set of polynomials, but these latter must be adapted to the structure of the sets $E_n$, especially when they converge to a Cantor set.

In \cite{dolo} we have introduced a different approach, by factoring the polynomial $Z(z)$ as in eq. (\ref{meas3}). The price to pay is non--linearity of the set of equations (\ref{meas5}), but the technique is stable and the roots $\zeta_i$, when suitably rescaled, provide a convenient parameterization of $Z(\zeta;z)$.
Eqs. (\ref{meas5}) can be solved \cite{dolo} by the hybrid Powell method (routine {\sc HYBRJ} in Minpack \cite{powell,minpack}). This technique evaluates the Jacobian matrix of the system of equations, thereby requiring a storage of $N^2$ quantities, and a computational complexity of at least $N^3$ operations. Since in this paper we let $n$ grow, to approach a limit Cantor set, while $N = M^n$, this is a severe limitation.

Luckily, the structure of the system (\ref{meas5}), when applied to Cantor sets, allows for a different algorithm, that does {\em not} require the computation of the Jacobian.
Observe that one can deal with the integrals (\ref{meas5}) on the {\em gaps} between intervals by Gaussian summation, as done in Sect. \ref{sec-equil} for the integral over the intervals $E^i_n$. In fact, for any $i = 1,\ldots,N-1$, let ${\cal K}_i$ be such integral, and let us factor the function $Y(s)$ as $Y(s) = (s-\beta_i)(\alpha_{i+1}-s) \tilde{Y_i}(s)$, so that
 \begin{equation}
 \label{meas5c}
    {\cal K}_i := \int_{\beta_i}^{\alpha_{i+1}} \frac{Z(\zeta;s)}{\sqrt{|Y(s)|}} \; ds =
    \int_{\beta_i}^{\alpha_{i+1}} \frac{Z(\zeta;s)}{\sqrt{|\tilde{Y_i}(s)|}}
    \frac{ds}{\sqrt{(s-\beta_i)(\alpha_{i+1}-s)}}
    \simeq \frac{1}{G} \sum_{l=1}^G \frac{Z(\zeta;\vartheta^i_l)}{\sqrt{|\tilde{Y_i}(\vartheta^i_l)|}},
 \end{equation}
where $\vartheta^i_l$, $l=1,\ldots,G$, are the Gaussian points of the Chebyshev measure on $[\beta_i,\alpha_{i+1}]$. Putting in evidence the $i$-th variable $\zeta_i$ in the $i$-th equation (\ref{meas5c}), we obtain
 \begin{equation}
 \label{meas5e}
{\cal K}_i \simeq
    \frac{1}{G} \sum_{l=1}^G
    (\vartheta^i_l - \zeta_i)
    \frac{\prod_{j \neq i} (\vartheta^i_l - \zeta_j) }{\sqrt{|\tilde{Y_i}(\vartheta^i_l)|}}  .
 \end{equation}
It was found in \cite{dolo} (Sect. 6, Fig. 2) that the ``diagonal'' derivatives $\frac{\partial {\cal K}_i}{\partial \zeta_i}$ largely exceed in magnitude non--diagonal ones. 
Therefore, we can exploit the simple structure of the equations (\ref{meas5e}) to set up an iterative algorithm for the solution of the non--linear system ${\cal K}_i=0$, $i = 1,\ldots,N-1$. Assume an initial set of values $\{\zeta^0_i\}_{i=1}^{N-1}$. For any $i$, solve ${\cal K}_i=0$, eq. (\ref{meas5e}), for $\zeta_i$, in terms of $\zeta_l=\zeta^0_l$, $l \neq i$. Update the solution vector $\zeta$ and iterate until convergence, that can be gauged by the absolute value of ${\cal K}$ and by the precision of the solution vector $\zeta$ \cite{dolo}. We have observed a rapid convergence of the technique, even in its simplest form---of course, one might also refine it by more sophisticate search techniques in the direction indicated by eqs. (\ref{meas5e}) \cite{powell}, and by letting $G$ be a function of the integration interval. It follows from the above that each step of this technique requires the evaluation of $2N$ Gaussian sums. The computational complexity of this technique is then proportional to $N G$, times the number of iteration required for convergence, that increases mildly with $n$, {\em i.e} the logarithm of $N$. It is therefore largely inferior to the other methods mentioned above. In addition, it can be easily programmed on a parallel machine, splitting the Gaussian sums among different processors. The large $N$ computations in this paper have been performed following this approach.

\section{Appendix B: reconsidering Gragg and Harrod's}
Gragg and Harrod's algorithm, that we have discussed in Sect. \ref{sec-gragg},
is termed RKPW (Rutishauer, Kahan, Pal, Walker)\footnote{Refer to \cite{gragg} Sect. 3, page 328 for notations and equations that we use below, observing that $\alpha_j$ corresponds to our $a_{j-1}$ and $\beta_j$ to $b_j$. Also, to stick with \cite{gragg}, we do not follow in this Appendix the coherent usage of dummy indices used in the remainder of the paper.}.  RKPW computes a Jacobi matrix of rank $G$, that is in one--to--one correspondence with the discrete measure $\sigma^G = \sum_{n=1}^G \omega_n D_{\lambda_n}$, composed of $G$  atoms, through $6 \; G^2$ arithmetic operations \cite{gragg}.
Yet, our problem is not to compute the {\em full} Jacobi matrix corresponding with a set of $G$ atoms.
Rather, we want to compute a {\em finite} truncation of this matrix, yet obtained with an {\em arbitrarily large} number $G$.

This can be obtained by a minor modification of RKPW, certainly obvious to its authors. Given the discrete measure $\sigma^G$, the algorithm adds one by one the atomic measures $\omega_n D_{\lambda_n}$ to the Jacobi matrix of the previous $n-1$ atoms, call it $J^{n-1}$. In so doing, the rank of the computed Jacobi matrix $J^n$ increases by one at each step. It requires a one--line proof to demonstrate theoretically that one can stop the computation of these matrices at any fixed truncation of rank $\bar{n}$, without affecting the result: in fact, the $\bar{n}$-truncation of $J^n$ (when $n \geq \bar{n}$) can be seen as the Jacobi matrix of a discrete measure with $\bar{n}$ atoms. This also
proves that the modified RKPW is also an algorithm to add a finite number of atomic measures to the Jacobi matrix of {\em any} arbitrary measure.

Numerically, this modification can be implemented as follows. In the original algorithm, the two sequences $\{a_j\}$, for $j=0$ to $n-1$ and $\{b^2_j\}$, for $j=0$ to $n-2$ (corresponding to $J^{n-1}$, {\em i.e.} to the first $n-1$ atoms) are entered in input. The entry $b^2_0$ is used to store the integral of the input measure. Furthermore, five auxiliary vectors are defined: $\rho$, $\gamma$, $\tau$, $\pi$ and $\nu$. They are all initialized to zero, except for $\gamma_0=1$ and $\pi_0^2=\omega_n$, where $\omega_n$ is the weight of the atomic measure, of location $\lambda_n$, that needs to be added. Notice that the input equation $\alpha_{n+1}=\lambda_n$ in \cite{gragg} is superfluous and can be omitted. The algorithm then proceeds by iteration: the step labeled by $k$ is used to compute the auxiliary vectors at index $k$ from their values at $k-1$. In this step, the updated values of $a_{k-1}$ and $b^2_{k-1}$ are also produced.

It is therefore not required to run the iteration for values of $k$ larger than the desired size $\bar{n}$ of the required Jacobi matrix, that is computed exactly for $k=\bar{n}$. This immediately cuts down the computational complexity to $12 \;\bar{n}G$ arithmetic operations and the storage requirement to nine vectors of {\em fixed} size $\bar{n}$, since the atoms $w_n D_{x_n}$ can be generated when required and need not to be stored.

This is the first improvement that we can bring to the basic RKPW algorithm.
Next, recall that the total number of atoms can be very large in our problem. We can further reduce the physical computation time of the truncated Jacobi matrix by exploiting the structure of the algorithm in a second way. The set of RKPW recursions only link vector values of index $k$ and $k-1$, and it runs in a stable way from $k=1$ to $k=\bar{n}$. We can therefore part the set of indices, and the corresponding vectors, among the $P$ processors of a parallel cluster with distributed memory. The flow of operations can be briefly described as follows: processor number 0 holds the initial chunk of the input Jacobi matrix. It starts the computation with the input $\lambda_l$ and $\omega_l$ (position and weight of the $l$-th atom). When the recurrence relation reaches the last index in processor 0, that is, $k=n/P$, it transmits the values of $\lambda_l$, $\omega_l$ and of the five auxiliary vector entries $\rho_k$, $\gamma_k$, $\tau_k$, $\pi_k$ and $\nu_k$, at $k=\bar{n}/P$, to the second processor. This latter operates on the second chunk of the Jacobi matrix. While it performs the recurrence relations in its range of indices, processor 0 can do the same with the next atom $\omega_{l+1}D_{\lambda_{l+1}}$. Of course, the same procedure can be extended to the full set of $P$ processors, resulting in a complete parallelization of the algorithm. The number of arithmetical operation per processor then scales as $A\bar{n}G/P$, and so does the physical time of the computation.


\begin{thebibliography}{}

\bibitem{ach} N.I. Akhiezer, {\em The Classical Moment Problem},
Hafner, New York, NY. (1965).

\bibitem{akhie1}
N. I. Akhieser,
{\em General theory of the P. L. Tchebycheff polynomials.}
Translated from the 1945 Russian original by B. Bojanov and G. Nikolov,
{\em East J. Approx.}  {\bf 12},  211--259 (2006).

\bibitem{akhie2}
N. I. Akhieser,
{\em Orthogonal polynomials on several intervals},
{\em Dokl. Akad. Nauk SSSR} {\bf 134} 9--12 (1960) (Russian); translated as
{\em Soviet Math. Dokl.} {\bf 1} 989--992 (1960)

\bibitem{andriev}
V.V. Andrievskii, {\em Constructive function theory on sets of the complex plane through potential theory and geometric function theory},
{\em Surveys in Approximation Theory}, {\bf 2} 1-52 (2006)

\bibitem{andriev2}
V.V. Andrievskii,
{\em On the Green function for a complement of a finite number of real intervals}, {\em Constr. Approx.} {\bf 20} 565--583 (2004).

\bibitem{andriev3}
V.V. Andrievskii,
{\em The highest smoothness of the Green function implies the highest density
of a set}, {\em Ark. Mat.} {\bf 42} 217--238 (2004).

\bibitem{sasha1} A.I. Aptekarev, {\em Asymptotic properties of polynomials orthogonal on a system of contours, and periodic motion of Toda lattices}, {\em Math. USSR Sb.} {\bf 53} (1986) 233--260.

\bibitem{ldos}
J. Avron, B. Simon,
{\em Almost periodic Schr\"odinger operators. II. The integrated density of states},
{\em Duke Math. J.} {\bf 50} (1983) 369--39.

\bibitem{baker}
G.A. Baker, D. Bessis, P. Moussa,
{\em A family of almost periodic Schr\"odinger operators},
{\em Phys. A} {\bf 124} 61--77 (1984).



\bibitem{ran1} L. Baribeau, D. Brunet, T. Ransford and J. Rostand, {\em Iterated function systems, capacity and Green's functions}, {\em Comput. Methods Funct. Theory} {\bf 4} (2004) 47--58.

\bibitem{ba2}  M. F. Barnsley, {\em Fractals Everywhere}, Academic Press,
New York, NY (1988).

\bibitem{dem}  M. F. Barnsley and S. G. Demko,
{\em Iterated function systems and the global construction of
fractals}, {\em Proc. R. Soc. London} \textbf{A 399} (1985)
243--275.

\bibitem{recur}
M.F. Barnsley, J. Elton, D.P. Hardin,
{\em Recurrent iterated function systems.
Fractal approximation},
{\em Constr. Approx.} {\bf 5} 3--31 (1989).

\bibitem{mijeff}
M. F. Barnsley, J. S. Geronimo and A. N. Harrington, {\em Almost periodic Jacobi matrices associated with Julia sets for polynomials}, {\em Comm. Math. Phys.} {\bf 99} 303--317 (1985).

\bibitem{belli} J. Bellissard, D. Bessis, and P. Moussa, {\em Chaotic states of almost--periodic Schr\"odinger operators}, {\em Phys. Rev. Lett.} {\bf 49} 701--704 (1982).

\bibitem{testard}
J. Bellissard, R. Lima, D. Testard,
{\em Almost periodic Schr\"odinger operators}, in {\em Mathematics + Physics} {\bf 1}, 1--64, World Sci. Publishing, Singapore, (1985).

\bibitem{bessi} D. Bessis, J. Geronimo, and P. Moussa, {\em Mellin transforms associated with Julia sets and physical applications}, {\em J. Stat. Phys.} {\bf 34} 75--110 (1983).

\bibitem{claudemichela} C. Brezinski, M. Redivo Zaglia,
{\em Extrapolation Methods: Theory and Practice}, North Holland, Amsterdam (1991).

\bibitem{brol} H. Brolin, {\em Invariant sets under iterations of rational functions}, {\em Ark. Mat.} {\bf 6} 103-144 (1965).

\bibitem{chris} J. S. Christiansen, B. Simon, M. Zinchenko, {\em Finite gap Jacobi matrices, I. The isospectral torus}, {\em Constr. Approx.} {\bf 32}, 1--65 (2010).

\bibitem{chris2} J. S. Christiansen, B. Simon, M. Zinchenko, {\em
Finite gap Jacobi matrices, III. Beyond the Szeg\"o class},
{\em Constr. Approx.} {\bf 35} 259--272 (2012).

\bibitem{chris3} J. S. Christiansen, {\em
Szeg\"o's theorem on Parreau-Widom sets},
{\em Adv. Math.} {\bf 229} 1180--1204 (2012).

\bibitem{cycon}
H.L. Cycon, R.G. Froese, W. Kirsch, B. Simon,
{\em Schr\"odinger operators with application to quantum mechanics and global geometry},
Springer-Verlag, Berlin (1987).

\bibitem{fibo2}
D. Damanik, A. Gorodetski,
{\em Spectral and quantum dynamical properties of the weakly coupled Fibonacci Hamiltonian}, {\em Comm. Math. Phys.} {\bf 305} 221--277 (2011).

\bibitem{damani} D. Damanik, B. Simon, {\em Jost function and Jost solution for Jacobi matrices, I}, {\em Invent. Math.} \textbf{165}, 1--50 (2006).
\bibitem{souil} F. Delyon and B. Souillard, {\em The rotation number for finite difference operators and its properties}, {\em Comm. Math. Phys.} {\bf 89}, 415--426 (1983).

\bibitem{tom}
E. Diekema, T.H. Koornwinder, {\em Differentiation by integration using orthogonal polynomials, a survey},
{\em J. Approx. Theory} {\bf 164}, 637-–667 (2012).

\bibitem{driscoll}
T.A. Driscoll, L. N. Trefethen,
{\em Schwarz-Christoffel mapping},
Cambridge Monographs on Applied and Computational Mathematics {\bf 8}, Cambridge University Press, Cambridge (2002).

\bibitem{nick} M. Embree and L. N. Trefethen,
{\em Green's Functions for Multiply Connected Domains via Conformal Mapping},
{\em SIAM Review} {\bf 41} 745--761 (1999).

\bibitem{fat} P. Fatou, {\em Sur les \'equations fonctionelles}, {\em Bull. Soc. Math. France} {\bf 47} 161--271 (1919), {\bf 48} 33--94 (1920), {\bf 48} 208--314 (1920).

\bibitem{hans2} H.-J. Fischer,
{\em Recurrence Coefficients of Orthogonal Polynomials with
Respect to Some Self-Similar Singular Distributions}, {\em Z. Anal. Anwendungen} \textbf{14}, 141--155 (1995).

\bibitem{hans3}
H.-J. Fischer, {\em On generating orthogonal polynomials for discrete measures}, {\em Z. Anal. Anwendungen}   {\bf 17}, 183--
205 (1998).


\bibitem{minpack} B. S. Garbow, K. E. Hillstrom and
J. J. More, Argonne National Laboratory. MINPACK project. March 1980.

\bibitem{gautschi}
W. Gautschi,
{\em Orthogonal polynomials: computation and approximation.}
Numerical Mathematics and Scientific Computation. Oxford Science Publications. Oxford University Press, New York, (2004)

\bibitem{jeff1}
J. S. Geronimo, W. Van Assche, {\em Orthogonal polynomials with asymptotically periodic recurrence coefficients}, {\em J. Approx. Theory} {\bf 46},  251--283 (1986).

\bibitem{jeff2}
J. S. Geronimo, W. Van Assche, {\em Orthogonal polynomials on several intervals via a polynomial mapping}, {\em Trans. Amer. Math. Soc.} {\bf 308}, 559--581 (1988).

\bibitem{gero}
Ya. L. Geronimus,
{\em On certain asymptotic properties of polynomials}. (Russian)
{\em Mat. Sbornik N. S.} {\bf 23} 77--88 (1948).

\bibitem{golw}
G. H. Golub and J. H. Welsch, {\em Calculation of Gauss
quadrature rules},  {\em Math. Comp.} {\bf 23}, 221--230,
1969.

\bibitem{gragg} W. B. Gragg and W. J. Harrod,
{\em The numerically stable reconstruction of Jacobi matrices from spectral data}, {\em Numer. Math.} \textbf{44}, 317--335 (1984).

\bibitem{martin}
M. H. Gutknecht,
{\em Solving Theodorsen's integral equation for conformal maps with the fast Fourier transform and various nonlinear iterative methods},
{\em Numer. Math.} {\bf 36}, 405--429 (1980/81).

\bibitem{doug}
D.P. Hardin, and E. B. Saff,
{\em Discretizing manifolds via minimum energy points},
{\em Notices Amer. Math. Soc.} {\bf 51} 1186--1194 (2004).

\bibitem{stric0}
S. M. Heilman, P. Owrutsky, and R. S. Strichartz,
{\em Orthogonal Polynomials with Respect to Self-Similar
Measures}, {\em Experiment. Math.}  {\bf 20}, (2011) 238--259.
\bibitem{henrici} P. Henrici, {\em Applied and computational complex analysis} {\bf 3}, John Wiley (New York) (1986).
\bibitem{hut}  J. Hutchinson, {\em Fractals and self--similarity},
{\em Indiana J. Math.}\ \textbf{30}, (1981) 713--747.

\bibitem{stric}
P. Janardhan, D. Rosenblum and R. S. Strichartz,
{\em Numerical experiments in Fourier asymptotics of Cantor measures and
wavelets},
{\em Experiment. Math.}  {\bf 1}, 249--273 (1992).

\bibitem{johnson} R. A. Johnson, {\em A review of recent work on almost
periodic differential and difference operators}, {\em Acta Appl. Math.} {\bf 1}, 241--261 (1983).

\bibitem{russell}
R. A. Johnson,
{\em Cantor spectrum for the quasi-periodic Schr\"odinger equation},
{\em J. Differential Equations} {\bf 91}, 88--110 (1991).

\bibitem{palle}
P. E. T. Jorgensen, K. A. Kornelson and K. L. Shuman, {\em Iterated Function Systems, Moments, and Transformations of Infinite Matrices}, {\em Memoirs of the AMS} {\bf 213} (2011).

\bibitem{jul} G. Julia, {\em M\'emoire sur l'iteration des fonctions rationelles}, {\em J. Math. Ser. 7 (Paris)} {\bf 4}, 47--245 (1918).

\bibitem{last} Y. Last, {\em Exotic spectra: a review of Barry Simon's central contributions}, in {\em Spectral theory and mathematical physics: a Festschrift in honor of Barry Simon's 60th birthday},
{\em Proc. Sympos. Pure Math.}  {\bf 76}, 697--712 (2007),
Amer. Math. Soc., Providence, RI.

\bibitem{last2}
Y. Last,
{\em Quantum dynamics and decompositions of singular continuous spectra},
{\em J. Funct. Anal.} {\bf 142}, 406--445 (1996).

\bibitem{dirk9} D. Laurie,
{\em Accurate recovery of recursion coefficients from Gaussian quadrature formulae}, {\em J. Comp. Appl. Math.}  \textbf{112}, 165--180 (1999).

\bibitem{lubin}
D. Lubinsky,
{\em Asymptotics of orthogonal polynomials: some old, some new, some identities},
{\em Acta Appl. Math.} {\bf 61}, 207--256 (2000).

\bibitem{alphonse}
A. Magnus,
{\em Recurrence coefficients for orthogonal polynomials on connected and nonconnected sets}, in {\em Pad\'e approximation and its applications (Proc. Conf., Univ. Antwerp, Antwerp, 1979},
{\em Lecture Notes in Math.}  {\bf 765}, 150--171, Springer, Berlin, (1979).

\bibitem{alphonse2}
A. Magnus,
{\em Pad\'e approximation to functions with branch points,}
Lecture notes 1998-1999,
http://perso.uclouvain.be/alphonse.magnus/num3/m3xxx98.ps

\bibitem{cap}  {G. Mantica}, {\em A Stieltjes Technique for Computing
Jacobi Matrices Associated With Singular Measures}, {\em Constr.
Appr.} {\bf 12}, (1996) 509--530.

\bibitem{physd1} {G. Mantica},
{\em Quantum intermittency in almost periodic systems derived from
their spectral properties}, {\em Physica D} {\bf 103}, (1997) 576--589.

\bibitem{physd2} {G. Mantica},
{\em Wave propagation in almost-periodic structures}, {\em Physica D}
{\bf 109}, (1997) 113--127.

\bibitem{etna}  G. Mantica,
{\em Fourier-Bessel functions of singular continuous measures and their
many asymptotics}, {\em Electron. Trans. Numer. Anal.} (Electronic) {\bf 25}, (2006) 409--430.

\bibitem{poin1}
 G. Mantica, S. Vaienti,
{\em The asymptotic behaviour of the Fourier transform of
orthogonal polynomials I: Mellin transform techniques},  {\em Ann.
Henri Poincar\'e} {\bf 8}, (2007) 265--300.

\bibitem{poin2}  G. Mantica, D. Guzzetti,
{\em The asymptotic behaviour of the Fourier transform of
orthogonal polynomials II: Iterated Function Systems and Quantum Mechanics},
{\em Ann. Henri Poincar\'e} {\bf 8}, (2007) 301--336.

\bibitem{nalgo2}  G. Mantica, {\em Dynamical Systems and Numerical Analysis: the Study of Measures generated by Uncountable I.F.S},  {\em Num. Alg.}  \textbf{55}, (2010) 321--335.

\bibitem{intjns}
G. Mantica, {\em On the attractor of one-dimensional infinite iterated function systems}, {\em Int. J.  Appl. non-linear Sci.} {\bf 1}, 87-99 (2013).

\bibitem{dolo}
G. Mantica, {\em Computing the equilibrium measure of a system of intervals converging to a Cantor set}, {\em DRNA} (Electronic) {\bf 6}, 51--61 (2013).

\bibitem{arxiv}
G. Mantica, {\em Direct and inverse computation of Jacobi matrices of infinite IFS}, {\em Numerische Math.} {\bf 125}, (2013)  705--731.

\bibitem{maroni} P. Maroni, {\em Une charact\'erisation des polyn\^omes orthogonaux semiclassiques}, {\em C. R. Acad. Sci. Paris S\'er. I Math.}
    {\bf 301} (1985) 269--272.

\bibitem{papmor}
P. A. P. Moran, {\em Additive functions of intervals and Hausdorff measure, Proc. Camb. Phil. Soc.} {\bf 42},  (1946) 15--23.

\bibitem{nuttall}
J. Nuttall and S. R. Singh, {\em Orthogonal polynomials and Pad\'e approximants associated with a
system of arcs, J. Approx. Theory} {\bf 21}, (1977) 1-–42.

\bibitem{franz0}
F. Peherstorfer, {\em On Bernstein-Szeg\"o orthogonal polynomials on several intervals}, {\em SIAM J. Math.
Anal.} {\bf 21}, (1990) 461–-482.

\bibitem{franz} F. Peherstorfer, {\em Orthogonal and extremal polynomials on several intervals}, {\em J. Comput. Appl. Math.} {\bf 48} (1993) 187--205.

\bibitem{franzimrn} F. Peherstorfer, {\em Zeros of polynomials orthogonal on several intervals}, {\em Int. Math. Res. Not.} {\bf 7} (2003) 361--385.

\bibitem{franzinfty} F. Peherstorfer, P. Yuditskii,
{\em Asymptotic behavior of polynomials orthonormal on a homogeneous set},
{\em J. Anal. Math.} {\bf 89} 113--154 (2003).

\bibitem{powell}
M. J. D. Powell,
{\em On nonlinear optimization since 1959}, {\em The birth of numerical analysis}, 141--160, World Sci. Publ., Hackensack, NJ, (2010).

\bibitem{pomme}
Ch. Pommerenke, {\em Boundary Behaviour of Conformal Maps}, Springer-Verlag, Berlin/New York, (1992).

\bibitem{ran} T. Ransford and J. Rostand, {\em Computation of Capacity}, {\em Math. of Comp.} {\bf 76} (2007) 1499--1520.
\bibitem{ran0} T. Ransford, {\em Potential theory in the complex plane}, Cambridge University Press, Cambridge (1995).

\bibitem{reedsim}
M. Reed, B. Simon,
{\em Methods of modern mathematical physics. IV. Analysis of operators}. Academic Press, Oxford University Press, New York  (2004).

\bibitem{lothar}
L. Reichel, {\em On polynomial approximation in the complex plane with application to conformal mapping},
{\em Math. Comp.} {\bf 44} 425--433 (1985);
{\em A fast method for solving certain integral equations of the first kind with application to conformal mapping}, Special issue on numerical conformal mapping, {\em J. Comput. Appl. Math.} {\bf 14} 125--142 (1986).

\bibitem{ed} E. B. Saff, {\em Logarithmic potential theory with applications to approximation theory}, {\em Surveys in Approx. Theory} {\bf 5} (2010) 165--200.
\bibitem{shen}
J. Shen, G. Strang, A. J. Wathen, {\em The potential theory of several intervals and its applications}, {\em Appl. Math. Optim.} {\bf 44} 67--85 (2001).
\bibitem{stahl} H. Stahl, V. Totik, {\em General Orthogonal Polynomials}, Cambridge University Press, Cambridge (2010).
\bibitem{simonratio} B. Simon, {\em Ratio asymptotics and weak asymptotic measures for orthogonal polynomials on the real line}, {\em J. Approx. Theory} {\bf 126} (2004) 198--217.

\bibitem{str0} R. S. Strichartz, {\em Analysis on Fractals},
{\em Notices of the AMS}, {\bf 46}, number 10, 1199--1208 (1999).

\bibitem{str0b}
R. S. Strichartz, {\em Differential Equations on Fractals:
A Tutorial}, Princeton University Press, (2006).

\bibitem{fibo1}
A. S\"uto,
{\em Singular continuous spectrum on a Cantor set of zero Lebesgue measure for the Fibonacci Hamiltonian}, {\em J. Statist. Phys.} {\bf 56} 525--531 (1989).

\bibitem{thouless}
D. J. Thouless,
{\em A relation between the density of states and range of localization for one dimensional random systems},
{\em J. Phys. C: Solid State Phys.},  {\bf 5}, 77--81 (1972).

\bibitem{serguey}
S. Tcheremchantsev, {\em Dynamical analysis of Schr\"odinger operators with growing sparse potentials}, {\em Comm. Math. Phys.} {\bf 253} (2005) 221--252.

\bibitem{totik}
V. Totik, {\em Metric properties of harmonic measures}, {\em Mem. Am. Math. Soc.}
{\bf 184} n. 867, (2006)

\bibitem{vanas}
W. Van Assche,
{\em Asymptotics for orthogonal polynomials and three-term recurrences}, in {\em Orthogonal Polynomials; Theory and Practice}, NATO-ASI series C {\bf 294},  435--462 (1990).

\bibitem{widom} H. Widom, {\em Extremal polynomials associated with a system of curves in the complex plane}, {\em Adv. in Math.} {\bf 3} (1969) 127--232.

\end{thebibliography}
\end{document}